\documentclass[11pt, oneside]{article}

\usepackage{amsmath}
\usepackage{amssymb}
\usepackage{graphics}
\usepackage{graphicx}
\usepackage{theorem}
\usepackage{latexsym}
\usepackage{setspace}

\usepackage{amssymb}

\newtheorem{lemma}{Lemma}
\newtheorem{theorem}{Theorem}
\newtheorem{conj}{Conjecture}

\newtheorem{prop}{Proposition}
\newtheorem{definition}{Definition}

\begin{document}

\title{Estimated transversality and rational maps}
\author{Rosa Sena-Dias}
\maketitle
\begin{abstract}
In this paper, we address a question of Donaldson's on the best estimate that can be achieved for the transversality of an asymptotically holomorphic sequence of sections of increasing powers of a line bundle over an integral symplectic manifold. More specifically, we find an upper bound for the transversality of $n$ such sequences of sections over a $2n$-dimensional symplectic manifold. In the simplest case of $S^2$, we also relate the problem to a well known question in potential theory (namely, that of finding logarithmic equilibrium points), thus establishing an experimental lower bound for the transversality.
\end{abstract}

\section{Introduction}
In his work on symplectic Lefschetz pencils, Donaldson introduced the notion of estimated transversality for a sequence of sections of a bundle tensored with increasing powers of a line bundle. This, together with asymptotic holomorphicity, is the key ingredient in the construction of symplectic submanifolds. Despite its importance in the area, estimated transversality has remained a mysterious property. In this paper, one of our aims is to shed some light into this notion by studying it in the simplest possible case, namely that of $S^2$. We state some new results about high degree rational maps on the 2-sphere that can be seen as consequences of Donaldson's existence theorem for pencils, and explain how one might go about answering a question of Donaldson: what is the best estimate for transversality that can be obtained? We also show how the methods applied to $S^2$ can be further generalized to prove the following:
\begin{theorem}\label{final} 
Let $X$ be a symplectic manifold of dimension $2n$ with symplectic form $\omega$ such that $[\omega/2\pi]$ lies in $H^2(X,\mathbb{Z})$, and a compatible almost complex structure. Let $L\rightarrow X$ be a Hermitian line bundle whose Chern class is $[\omega/2\pi]$. There exits $\eta_0<1$ such that, if we have $n$ asymptotically holomorphic sequences of sections $s_0,\cdots, s_n$ of $L^k$ satisfying
\(
\eta \leq ||s_0||^2+\cdots+||s_n||^2 \leq 1,
\)
then $\eta<\eta_0$.
\end{theorem}
Donaldson's results (see \cite{do1}, \cite{do2}) ensure that such sequences exist for some $\eta$ and any choice of complex structure. Recall that a sequence of sections $\{s_k\}$ of $L^k$ is said to be asymptotically holomorphic if its $\bar{\partial}$ is bounded independently of $k$. 

It is not hard to see that Donaldson's transversality theorem for the symplectic manifold $(S^2, \omega_{FS})$ implies the following interesting result:
\begin{prop}\label{principal}
There is $0<\eta \leq 1$, such that, for each $k$ large enough, there exists a pair of homogeneous polynomials, $(p_k,q_k)$, of degree $k$, in two complex variables, defining a function from $\mathbb{C}^2$ to $\mathbb{C}^2$, that takes $S^3\subset \mathbb{C}^2$ into an annulus of outer radius 1 and inner radius $\sqrt{\eta}$.
\end{prop}
This proposition is a consequence of a more general result for complex K\"ahler manifolds that comes from applying Donaldson's techniques to the complex setting (using the techniques appearing in \cite{do1} for the K\"ahler setting; this is an exercise which we carry out for the sake of completeness). The special case of $S^2$ is treated in more detail in section \ref{capitulo_S2}. To get a feeling for how strong the above result is, let us try to take $p_k(\mathbf{z},\mathbf{w})=\mathbf{z}^k$ and $q_k(\mathbf{z},\mathbf{w})=\mathbf{w}^k$. Then, the image of $S^3$ by each of the maps $(p_k,q_k)$ is contained in an annulus of outer radius 1 but whose inner radius is $1/\sqrt{2^{k-1}}$. 

In fact, a sequence of pairs of homogeneous polynomials satisfying Proposition \ref{principal} is very special. We prove:
\begin{theorem}\label{ud}
Let $(p_k,q_k)$ be a sequence of pairs of homogeneous polynomials as above. The map $p_k/q_k$, thought of as a degree $k$ map of $\mathbb{CP}^1$ to itself, has asymptotically uniformly distributed fibers, in the sense that, if $x^k_i$ denote the points in one fiber for $i=1, \cdots k$, counted with multiplicity, and $f$ is a $\mathcal{C}^2$ function on $S^2$, then
\begin{equation}
\left| \frac{1}{k}\sum_{i=1}^k f(x^k_i)-\frac{1}{|S^2|}\int_{S^2} f \right|\leq \frac{C||\triangle f||_{\infty}}{k}\label{unif_distribuido},
\end{equation}
and, in particular, tends to zero. 
\end{theorem}
A similar result holds for the branch points of $p_k/q_k$. The problem of distributing points on $S^2$ is an old and important problem with many applications. It has been addressed by several branches of mathematics, for example in potential theory (see \cite{rsz}) and in arithmetic number theory (see \cite{s}). This result is sharper than results found through potential theory or arithmetic number theory methods (although in this last case the bound for the expression in inequality (\ref{unif_distribuido}) is for functions $f$ in $L^2$ and not simply in $\mathcal{C}^2$). Indeed the bound (\ref{unif_distribuido}) is optimal in the sense that we cannot expect to get a better asymptotic bound in $k$ for $\mathcal{C}^2$ functions using second derivatives. 
The statement of the above Proposition \ref{principal}, has nothing to do with symplectic geometry, it is simply a statement about rational maps. We try to look at it without using the techniques of \cite{do2} and give an explicit construction of polynomials which are experimentally seen to satisfy the required property. This involves the choice of two sets of $k$ points, the zeroes of $p_k$ and $q_k$ on $S^2$. We will choose two sets of asymptotically uniformly distributed points which are a slight modification of the so called generalized spiral points. Generalized spiral points come from trying to solve a problem in potential theory, that of distributing a big number of charges on the 2 sphere subject to a logarithmic potential. They are described in \cite{rsz} as a good approximation of the actual solution to this problem. Even though the problem itself remains unsolved, some things are known about the optimal distribution. We will discuss the relations between this problem and our own. One of the upshots of this explicit construction is that it allows to experimentally determine a lower bound for the constant $\eta$ appearing above. Together with the upper bound coming from Theorem \ref{final} (that can be made totally explicit for $S^2$), this gives a partial answer to Donaldson's question.

A brief outline of this paper is the following: In section \ref {capitulo_transversalidade}, we give a review of the results in \cite{do1} and \cite{do2} and explain how they work  in the complex setting. In section \ref{capitulo_S2}, we explain how these results imply Proposition \ref{principal} and we also calculate an explicit upper bound for $\eta$ appearing in that Proposition, this bound turns out to be of the order $e^{-10^{36}}$! We also prove Theorem \ref{ud}. Section 4 gives a construction of a sequence of polynomials which are experimentally shown to satisfy the condition in Proposition \ref{principal}, as well as some steps towards the proof of the fact that they indeed satisfy the required condition. We also describe some relations between estimated transversality and the logarithm equilibrium problem in potential theory. Section 5 generalizes the method described in section 3 to find an upper bound for $\eta$ in the general case, thus proving Theorem \ref{final}. We end with some conjectures on how to find polynomials satisfying Proposition \ref{principal} via PDE theory.

{\bf Acknowledgments} I am very grateful to my advisor, Peter Kronheimer, for his guidance and support throughout my Ph.D. and also for the many helpful discussions on the subject of this paper. I would also like to thank Pedro for being an inexhaustible source of joy and peace to me.

\section{Background}\label{capitulo_transversalidade}
The notion of linear system is one of importance in complex geometry. For example, it is often using a linear system that one is able to realize a complex manifold (when it satisfies certain constrains) as a submanifold of $\mathbb{CP}^N$.  This is the content of Kodaira's embedding theorem. More trivially, a generic linear system of dimension 0 (whose existence is easy to establish) gives rise to a divisor, i.e., a complex submanifold. Moving one step up, a generic one dimensional linear system gives rise to a Lefschetz pencil, i.e., a holomorphic map $X\rightarrow \mathbb{CP}^1$ (defined away from a codimension $4$ subvariety) with the simplest possible singularities. It is well know, that every symplectic manifold has an almost-complex structure, which is the same as saying, that the symplectic category generalizes the K\"ahler category. The natural question is then: is it possible to generalize the notion of linear system to this new setting, and use it to study symplectic manifolds, as it was used to study complex manifolds? Even in the simplest case of linear systems of dimension $0$, this poses problems. The most important difficulty with which one is faced is the non-existence of holomorphic sections of complex bundles over symplectic manifolds (except in the integrable case where the manifold is actually complex). In his paper \cite{do1}, Donaldson, resolves this issue by substituting the holomorphic condition by what is called asymptotic holomorphicity. One looks for sections of an increasing power of a line bundle which have a bounded $\bar\partial$. In the holomorphic setting, a complex submanifold is obtained as the zero set of a transverse, holomorphic section. The holomorphicity condition becomes asymptotic holomorphicity, how about the transverse condition? Even though this condition could easily be translated into the symplectic picture, it is no longer strong enough. It needs to become "estimated transversality", that is, transverse with a good estimate independent of the (increasing) degree of the bundle. Using this notion in a very key manner, Donaldson establishes the existence of sections of bundles whose zero set is symplectic, therefore proving an important existence theorem for symplectic submanifolds. In \cite{do2}, he goes one step further and proves the existence of the analogue of pencils. These symplectic Lefschetz pencils completely characterize symplectic manifolds: very roughly, a manifold is symplectic exactly when it can be seen (after blow up) as a bundle over $\mathbb{CP}^1$, with symplectic fibers, some of which have simple singularities. Surprisingly, the techniques used to prove this, and in particular the notion of estimated transversality, can be used to prove new theorems in the complex setting. Even for $S^2$, one can prove new and unexpected results on rational maps of high degree.

\subsection{Symplectic Lefschetz pencils and estimated transversality}
 A natural question to ask in symplectic geometry is: Does every symplectic manifold have symplectic submanifolds? There have, so far, been two approaches to this question:
\begin{enumerate}
\item The oldest one, by Gromov, is, given an almost-complex structure $J$ on $X$, to look for 2-dimensional submanifolds which are $J$ invariant. More precisely, $Y$ a dimension 2 submanifold of $X$ is a $J$ holomorphic curve if, for all $y\in Y$, $J(T_yY)=T_yY$. Every holomorphic curve is symplectic.
 \item The more recent one, the one we are going to be concerned with here, by Donaldson, is, assuming that $[\omega/2\pi]$ is an integral class and that $L\rightarrow X$ is a complex line bundle with $c_1(L)=[\omega/2\pi]$, to look for sections of the bundle $L^k\rightarrow X$ with $|\bar{\partial}s|<|\partial s|$ (here $\partial$ and $\bar\partial$ are with respect to some almost-complex structure). Their zero set will be a codimension 2 symplectic submanifold of $X$. 
\end{enumerate}
To be more precise, in his paper \cite{do1}, Donaldson proves the existence of such sections:
\begin{theorem}[Donaldson, \cite{do1}]\label{existe_secc}
 Let $X$ be a manifold with a symplectic form $\omega$, such that $[\omega/2\pi]\in H^2(X,\mathbb{Z})$ and let $L\rightarrow X$ be a complex Hermitian line bundle with a connection form, whose curvature is $i\omega$. For sufficiently large $k$, there is a sequence, ${s_k}$ of sections of $L^k$ such that:
 \begin{enumerate}
 \item $|s_k|$ is bounded by 1, $|\nabla s_k|\leq C\sqrt{k}$ and $|\nabla \nabla s_k|\leq Ck$ where C is independent of $k$, 
 \item \label{asym_holo} $|\bar\partial s_k|$ is bounded by some constant $C$, independent of $k$,
 \item  there is a constant $\eta$, independent of $k$ such that  $|s_k|\leq \eta \implies |\partial s_k|\geq \eta\sqrt{k}$.
 \end{enumerate}
 \end{theorem} 
When $k$ is large enough, $C< \eta\sqrt{k}$ along the zero set of $s_k$, therefore $|\bar\partial s_k|<|\partial s_k|$ which implies that $s_k^{-1}(0)$ is a submanifold of $X$ and that it is symplectic. Condition 3 in Theorem \ref{existe_secc} plays an extremely important role in the story. A sequence satisfying it is said to be $\eta$ transverse to zero (or simply $\eta$ transverse). In fact, we can generalize this notion further:
\begin{definition}
Let $X$ be a manifold with a metric, $L\rightarrow X$ a Hermitian complex line bundle, $E\rightarrow X$ a Hermitian complex vector bundle and $\{\tau_k\}$ a sequence of sections of $E\otimes L^k$. Let $\eta$ be a positive number. Then, we say that $\{\tau_k\}$ is $\eta$ transverse to zero if
\begin{displaymath}
|\tau_k|\leq \eta \implies \langle[\partial \tau_k]^*\mathbf{v},[\partial \tau_k]^*\mathbf{v}\rangle\geq \eta^2{k}|\mathbf{v}|^2,
\end{displaymath}
for all $\mathbf{v}$ section of $E\otimes L^k$.
\end{definition}
Note that, if $\partial \tau_k$ is not surjective, this will not be possible, since for some $\mathbf{v}\ne 0$, $[\partial \tau_k]^*\mathbf{v}=0$. The above definition is the same as asking that $|\tau_k|\leq \eta$, needs to imply that $\partial \tau_k$ is surjective and has a pointwise right inverse, whose norm is smaller than $\eta^{-1}k^{-1/2}$. Condition 2 in Theorem \ref{existe_secc} is referred to as asymptotic holomorphicity.

Generalizing this result further, Auroux in \cite{a1}, proves Theorem \ref{existe_secc} with $L$ replaced  by $L\otimes E$, where $E$ is a complex vector bundle of any rank. Theorem \ref{existe_secc} is an existence theorem for the symplectic analogs of linear systems of dimension $0$. What about pencils, i.e., linear systems of dimension $1$? The first thing to do is to generalize the notion to this new setting. 
\begin{definition}
Let $X$ be a manifold with a symplectic form. Then, a map $F$, defined on $X$ minus a codimension 4 manifold, is a symplectic Lefschetz pencil if, for every point $p\in X$, one of the following conditions is satisfied:
\begin{itemize}
\item F is a submersion at $p$,
\item F is not defined at $p$, in which case, there are compatible complex coordinates $z_1,\cdots z_n$, centered at $p$, with $F={z_1}/{z_2}$, 
\item F is defined at $p$, but it is not a submersion at $p$, there are compatible complex coordinates $z_1,\cdots z_n$, centered at $p$, with $F=z_1^2+\cdots z_n^2$.
\end{itemize} 
\end{definition}
Here, "compatible complex coordinates" simply means a map, defined locally around $p$ to $\mathbb{C}^n$, such that, the pullback of the standard symplectic form on $\mathbb{C}^n$ is $\omega$ at the origin. Donaldson in \cite{do2}, proves the following:
\begin{theorem}[Donaldson, \cite{do2}]\label{exist_lapiz}
 Let $X$ be a manifold with a symplectic form $\omega$, such that $[\omega/2\pi]\in H^2(X,\mathbb{Z})$ and let $L\rightarrow X$ be a complex Hermitian line bundle with a connection form, whose curvature is $i\omega$. For sufficiently large $k$, there is a sequence of pairs of sections $(s_{0},s_{1})$ of $L^k$ such that:
\begin{enumerate}
\item $|s_0|^2+|s_1|^2\leq 1$, $|\nabla s_i|\leq C\sqrt{k}$ and $|\nabla\nabla s_i|\leq C{k}$, $i=0,1$,
\item $|\bar\partial s_0|$ and $|\bar\partial s_1|$ are bounded by a constant independent of $k$,
\item $s_0$ is $\eta$ transverse to zero, for some $\eta$ independent of $k$,
\item $(s_0,s_1)$ is $\eta$ transverse to zero, 
\item $\partial({s_1}/{s_0})$ is $\eta$ transverse to zero, away from the zero locus of $(s_0,s_1)$\label{derivada_transversal}.
\end{enumerate}
\end{theorem}
As explained in \cite{do2}, for large $k$, after perturbing the sections $s_1$ slightly, the map ${s_1}/{s_0}$ will give rise to a Lefschetz pencil. 

\subsection{Estimated transversality in the complex setting}\label{seccao_estimativa_complexa}
Since any K\"ahler manifold is also symplectic, one could apply Theorems \ref{existe_secc} and \ref{exist_lapiz} to the K\"ahler case. At first sight, these theorems seem not to generalize the existence theorems for complex submanifolds and holomorphic Lefschetz pencils, since they do not produce holomorphic sections of bundles. But Donaldson, in \cite{do1}, proves that in the K\"ahler case, the asymptotically holomorphic condition in Theorem \ref{existe_secc} can be strengthened to holomorphic. The result then becomes a new theorem for K\"ahler manifolds: 
\begin{theorem}[Donaldson, \cite{do1}]\label{existe_uma_holomorfa}
 Let $X$ be a K\"ahler manifold with integral cohomology. Let $\omega$ be the symplectic form on $X$ and let $L\rightarrow X$ be a complex Hermitian line bundle with a connection form, whose curvature is $i\omega$. For sufficiently large $k$, there is a sequence of holomorphic sections of $L^k$, $\{s_k\}$, such that:
 \begin{itemize}
 \item $|s_k|$ is bounded by 1,
 \item  \label{transversal} there is a constant $\eta$, independent of $k$, such that  $|s_k|\leq \eta \implies |\partial s_k|\geq \eta\sqrt{k}$.
 \end{itemize}
 \end{theorem} 
It is then natural to ask: what does Theorem \ref{exist_lapiz} say for a K\"ahler manifold?
\begin{theorem}\label{exist_lapiz_holo}
Let $X$ be a K\"ahler manifold with symplectic form $\omega$, such that $[\omega/2\pi]$ is in $H^2(X,\mathbb{Z})$ and a complex Hermitian line bundle $L\rightarrow X$ with a connection form, whose curvature is $i\omega$. For sufficiently large $k$, there is a sequence of holomorphic sections of $L^k$, $(s_0,s_1)$, such that:
\begin{enumerate}
\item $|s_0|^2+|s_1|^2\leq 1$, $|\nabla s_i|\leq C\sqrt{k}$ and $|\nabla\nabla s_i|\leq C{k}$, $i=0,1$,\label{sup}
\item $s_0$ is $\eta$ transverse to zero, for some $\eta$ independent of $k$,\label{uma_transversal}
\item $(s_0,s_1)$ is $\eta$ transverse to zero, \label{par_transversal}
\item $\partial({s_1}/{s_0})$ is $\eta$ transverse to zero, away from the zero locus of $(s_0,s_1)$.
\end{enumerate}
\end{theorem}
The proof of Theorem \ref{exist_lapiz_holo} is a combination of elements in the proof of Theorem \ref{exist_lapiz} and elements in the proof of Theorem \ref{existe_uma_holomorfa}. The first step is to get $s_0$, a section of $L^k$ for big $k$, which is holomorphic and $\eta$ transverse to zero, for some $\eta$ independent of $k$, by using Theorem \ref{existe_uma_holomorfa}. Next, one needs to build a sequence of pairs of sections of $L^k$, that are $\eta$ transverse to zero (note that this pair of sequences can be built close to $(s_0,0)$, so that its first term will be $\eta/2$ transverse). Instead of using the asymptotically holomorphic pair built in the proof of Theorem \ref{exist_lapiz} and perturbing it to make it holomorphic, as in the proof of Theorem \ref{existe_uma_holomorfa}, one can use the methods of \cite{do2} for building pencils, together with the following existence lemma proved in \cite{do1}:
\begin{lemma}[Donaldson,\cite{do1}]\label{exits_sigma_p}
There are constants $a$, $b$ and $c$ such that, given any $p\in X$ and any $k$ large, there is a holomorphic section $\sigma_p$ of $L^k$, satisfying the following estimates:
\begin{itemize}
\item $e^{-bkd^2(p,q)} \leq  | \sigma_p(q)| \leq  e^{-akd^2(p,q)}$ if $d(p,q)\leq ck^{-1/3}$, 
\item $|  \sigma_p(q)| \leq  e^{-ak^{1/3}}$ if $d(p,q)\geq ck^{-1/3}$.
\end{itemize}
The section $\sigma_p$ also satisfies
\begin{itemize}
\item $|\partial \sigma_p(q)| \leq \sqrt{k} e^{-akd^2(p,q)}$ if $d(p,q)\leq ck^{-1/3}$, \label{desig}
\item $|\partial  \sigma_p(q)| \leq  \sqrt{k} e^{-ak^{1/3}}$ if $d(p,q)\geq ck^{-1/3}$.
\end{itemize}
\end{lemma}
To proceed, the strategy is the usual one. Use a covering lemma to cover $X$ with colored neighborhoods of points where the above lemma applies (two neighborhoods with the same color are well separated). Now start with any pair of holomorphic sections, $(s_0,0)$ for example, and modify it, over all the balls of the first color, by adding a section of the form $(w_1 \sigma_{p_i}, w_2 \sigma_{p_i})$ for each ball (of color $1$). The vector $w=(w_1,w_2)$ comes from applying a transversality lemma to the representation  of the pair with respect to the trivialization $\sigma_{p_i}$. The lemma is:
\begin{lemma}[Donaldson,\cite{do2}]\label{trans}
Given a map $f:B^m(11/10) \rightarrow  B^n(1)$, which is holomorphic and any small positive number $\delta$, there is a $w\in \mathbb{C}^n$ with norm smaller than $\delta$, such that $f+w$ is $\eta={\delta}/{\log^p{\delta}}$ transverse to zero.
\end{lemma}
In this way, we achieve both holomorphicity and transversality over each ball in the first color. Note that the balls of a given color do not interfere with each other since they are well separated (here we use the inequalities in Lemma \ref{exits_sigma_p}). Next, we apply the same method to the balls of the second color, thus achieving transversality there. But we need to make sure that we do not spoil the transversality achieved over the balls with the first color. This is a consequence of the inequalities in Lemma \ref{exits_sigma_p}. We can keep up this process until we have gone trough all the colors. Just as in the proof of Theorem \ref{exist_lapiz}, (see \cite{do2}), the method produces an $\eta$ transverse sequence of sections, for some $\eta$ independent of $k$. These sequences are actually holomorphic, simply because they are a sum of holomorphic sections. Starting from the $\eta$ transverse pair, we slightly modify the second term $s_1$, to obtain transversality for $\partial(s_1/s_0)$, just as in \cite{do2}, but using holomorphic reference peak sections. We make use of the following lemma: 	
\begin{lemma}\label{trans_pi}
Let $\delta$ be a positive number smaller than $1/2$ and $p$ a point in $X$. For sequences $s_0$ and $s_1$ as above, large $k$ and small $r$, there is a $\pi\in \mathbb{C}^n$ (not depending on $k$) with norm smaller than $\delta$, such that, 
\begin{displaymath}
\partial \left(\frac{s_1+\sigma_{p,\pi}}{s_0}\right)
\end{displaymath}
is ${\delta}/{\log^p({1/\delta})}$ transverse to zero, over the ball of center $p$ and radius $rk^{-1/2}$. Here, 
\begin{displaymath}
\sigma_{p,\pi}=\sum_\alpha \pi^\alpha \sigma^\alpha_p
\end{displaymath}
and $\sigma^\alpha_{p}$ comes from Lemma \ref{exits_sigma_p}.
\end{lemma} 
This can be proved exactly as the corresponding statement in \cite{do2}. Now we can proceed as in \cite{do2}. The final result will be holomorphic, because the $\sigma_{p,\pi}$ are. Note that a modification of the above construction for the pair can be used to prove the following:
\begin{prop}\label{existem_n} 
Let $X$ be a K\"ahler manifold, of complex dimension $n$, with symplectic form $\omega$, such that $[\omega/2\pi]$ lies in $H^2(X,\mathbb{Z})$. Let $L\rightarrow X$ be a complex Hermitian line bundle $L\rightarrow X$ with a connection form, whose curvature is $i\omega$. There exits $0<\eta\leq 1$ and $n$ sequences of holomorphic sections $s_0,\cdots, s_n$ of $L^k$ satisfying:
\(
\eta \leq ||s_0||^2+\cdots+||s_n||^2 \leq 1.
\)
\end{prop}

\section{Estimated transversality for rational maps on $S^2$}\label{capitulo_S2}

\subsection{The problem on $S^2$}
In \cite{do2}, Donaldson asks the question: Given a symplectic manifold, what is the best $\eta$ for which we can find $(s_0,s_1) \in \Gamma(L^k)$ satisfying the conditions of Theorem \ref{exist_lapiz}? (We need to normalize the pairs to have $L^\infty$ norm 1 for this question to make sense.) We will address this question for $S^2$ with the Fubini-Study metric. Now, $S^2$ with the Fubini-Study form is K\"ahler, so we can apply to it Theorem \ref{exist_lapiz_holo}. In fact, the existence of a pair of holomorphic sections of $\mathcal{O}(k)$, satisfying only conditions \ref{sup} and \ref{par_transversal} in the theorem, seems to already give an interesting result. This last condition becomes somewhat simpler in the context of 2-dimensional manifolds, for there can be no surjective maps from $\mathbb{C}^2$ to $\mathbb{C}$. Namely, it becomes $|s_0|^2+|s_1|^2\geq \eta$ so that, when $|s_0|^2\leq \eta/2$, $|s_1|^2\geq \eta/2$. Also, holomorphic sections of $\mathcal{O}(k)$ are easy to characterize, they are simply homogeneous polynomials of degree $k$ in two complex variables. In this way, we prove Proposition \ref{principal}. A way to find a lower bound for the best $\eta$ appearing in that proposition is then to explicitly determine a sequence of pairs of homogeneous polynomials of degree $k$, $(p_k,q_k)$, such that the number 
\begin{displaymath}
\frac{\max (||p_k ||^2+||q_k||^2)}{\min (||p_k||^2+||q_k||^2)}
\end{displaymath}
is bounded independently of $k$. Where $|| p_k||$ stands for the norm of $p_k$ as a section of $\mathcal{O}(k)$, i.e., letting $[\mathbf{z}:\mathbf{w}]$ be the homogeneous coordinates in $S^2=\mathbb{C}\mathbb{P}^1$, 
\begin{displaymath}
||p_k||[\mathbf{z}:\mathbf{w}]=\frac{|p_k|(\mathbf{z},\mathbf{w})}{(|\mathbf{z}|^2+|\mathbf{w}|^2)^{k/2}}.
\end{displaymath}
The inverse of this bound will provide the lower bound we are looking for. There is a chart on $S^2$ minus the south pole, obtained by stereographic projection trough the south pole. It is centered at the north pole and identifies $S^2$ minus the south pole with $\mathbb{C}$. In homogenous coordinates on $S^2$, it is simply given by $z=\mathbf{z}/\mathbf{w}$. Over $S^2$ minus the south pole, $\mathcal{O}(k)$ admits a trivialization, i.e., a global section which we denote by $\mathbf{w}^k$. In fact, when we identify sections of $\mathcal{O}(k)$ with homogeneous polynomials of degree $k$, this section corresponds precisely to the polynomial $\mathbf{w}^k$. It is actually defined over all of $S^2$, but vanishes at the south pole. Its norm in the $z$ coordinate is simply,  
\begin{displaymath}
\frac{1}{(1+|z|^2)^{k/2}}.
\end{displaymath}
Often, we do not distinguish between a homogeneous polynomial and its representation in this trivialization, which is simply a polynomial of degree $k$, in one complex variable $z$.
 
\subsection{An upper bound for $\eta$ on $S^2$}\label{eta_para_S2}
In this section, we show how to find an upper bound for the best $\eta$ for which there is a sequence of pairs of polynomials, $(p_k,q_k)$, mapping $S^3$ into the $(\eta,1)$ annulus in $\mathbb{C}^2$. 
\begin{prop}\label{upper_bound}
If $p_k$ and $q_k$ are two homogeneous polynomials of degree $k$ in two complex variables, and if, for $k$ large, $(p_k,q_k)$ maps $S^3$ into the annulus of outer radius 1 and inner radius $\eta$ (for some $eta$ independent of $k$) then, $\eta \leq\eta_0 < 1$, where $\eta_0$ is going to be made explicit ahead.
\end{prop}
{\bf proof}:
Choose complex coordinates on $S^2$, by using stereographic projection trough the north pole for example. Consider $p_k$ and $q_k$ as polynomials in the chosen complex coordinate. We set
\begin{displaymath}
F_k=\frac{p_k}{q_k},
\end{displaymath}
\begin{displaymath}
u_k=\frac{1}{2}\log \left( \frac{F^*_k\omega_{FS}}{\omega_{FS}}\right),
\end{displaymath}
(note that $u$ has log type singularities) and 
\begin{displaymath}
f_k=\log(||p_k||^2+||q_k||^2).
\end{displaymath}
We have
\begin{displaymath}
\frac{{F_k}^*\omega_{FS}}{\omega_{FS}}=\frac{|p'_k q_k-p_kq'_k|^2(1+|z|^2)^2}{(|p_k|^2+|q_k|^2)^2},
\end{displaymath}
which can be written as
\begin{displaymath}
\frac{F_k^*\omega_{FS}}{\omega_{FS}}=\frac{|p_k'q_k-p_kq_k'|^2}{(1+|z|^2)^{2k-2}}\frac{(1+|z|^2)^{2k}}{(|p_k|^2+|q_k|^2)^2}
\end{displaymath}
and therefore
\begin{displaymath}
\frac{F_k^*\omega_{FS}}{\omega_{FS}}=\frac{||s_k||^2}{(||p_k||^2+||q_k||^2)^2},
\end{displaymath}
where $s_k$ is a section of $\mathcal{O}(2k-2)$, so that
\begin{displaymath}
u_k=\frac{1}{2}\log(||s_k||^2)-f_k.
\end{displaymath}
Now, we can use the fact that $\triangle \log ||s_k||^2=-(2k-2)$, away from the zeroes of $s_k$, to conclude that
\begin{displaymath}
\triangle u=-(k-1)-\triangle f_k,
\end{displaymath}
(we dropped the $k$ dependence on $u$). Now, $e^{2u}$ is the conformal factor of the metric defined by $F_k^*\omega_{FS}$. Note that $e^{2u}$ actually has zeroes. Expressing the curvature of the pullback metric in terms of this factor and noting that this curvature is one away from the branch points of $F_k$, one concludes that $u$ satisfies the following PDE, away from the branch points of $F_k$,
\begin{displaymath}
\triangle u+e^{2u}=1,
\end{displaymath}
that is, 
\begin{displaymath}
\triangle u+ \frac{F_k^*\omega_{FS}}{\omega_{FS}}=1.
\end{displaymath}
Putting this together with $\triangle u=-(k-1)-\triangle f_k $, we conclude
\begin{equation}\label{PDE}
\frac{1}{k}\frac{F_k^*\omega_{FS}}{\omega_{FS}}=1+\frac{1}{k}\triangle f_k,
\end{equation}
which holds at all points on $S^2$. Next, we show that the polynomials $p_k$, $q_k$, after rescaling, satisfy bounds independent of $k$, as do all of their derivatives. We start by recalling that
\begin{equation}\label{entre_eta_e_1}
\eta \leq \frac{|p_k|^2(z)+|q_k|^2(z)}{(1+|z|^2)^k}\leq 1.
\end{equation}
Now consider the rescaled polynomials
\begin{displaymath}
\tilde{p_k}(z)=p_k(\frac{z}{\sqrt{z}}), \,\,\,\,\tilde{q_k}(z)=q_k(\frac{z}{\sqrt{z}}),
\end{displaymath}
as a consequence of inequality \ref{entre_eta_e_1}, these satisfy 
\begin{displaymath}
\eta \left(1+\frac{|z|^2}{k} \right)^k \leq |\tilde{p}_k|^2+|\tilde{q}_k|^2\leq \left (1+\frac{|z|^2}{k}\right )^k,
\end{displaymath}
so that, $|\tilde{p}_k|$ and $|\tilde{q}_k|$, when restricted to the disk of center $0$ and radius $2c$, $D_{2c}$, admit an upper bound $e^{2c^2}$ (also valid on the closure of the disc). Now both polynomials are holomorphic functions on  $D_{2c}$. The Cauchy formula gives an upper bound for the norm of their derivatives on $D_c$, namely $2e^{2c^2}$. In fact, we even get an upper bound for the norm of higher order derivatives. If $l$ denotes the order of the derivative, 
\begin{displaymath}
|\tilde{p}_k^l|,|\tilde{q}_k^l| \leq \frac{2l!e^{2c^2}}{c^{l-1}},
\end{displaymath}
and
\begin{displaymath}
|\tilde{p}_k|,|\tilde{p}_k| \leq e^{\frac{c^2}{2}}.
\end{displaymath}
This then proves that $\tilde{p}_k$ and $\tilde{q}_k$ have convergent subsequences. Let us call the limits $p$ and $q$ respectively. These satisfy
\begin{displaymath}
\eta \leq |p|^2+|q|^2 \leq e^{c^2}.
\end{displaymath}
Also, the derivatives of $p$ and $q$ satisfy the same bounds as those of $\tilde{p}_k$ and $\tilde{q}_k$. From the inequality $\eta \leq |p|^2+|q|^2 $, not both $p$ and $q$ can vanish at zero. In fact, either $|p(0)|$ or $|q(0)|$ is greater than $\sqrt{\eta}/\sqrt{2}$. Let's say that $q$ satisfies this bound. Then, if $c$ is such that 
\begin{displaymath}
2e^{2c^2}c \leq \frac{\sqrt{\eta}}{2\sqrt{2}},
\end{displaymath}
in $D_c$,
\begin{displaymath}
|q|\geq  \frac{\sqrt{\eta}}{2\sqrt{2}}.
\end{displaymath}
This comes from the lower bound on $|q(0)|$ and the upper bound on $|q'|$. From now on, we assume that $c$ is such that it satisfies this. Let $G_k$ be the quotient of $\tilde{p}_k$ and $\tilde{q}_k$ (that is $F_k$ rescaled). On $D_c$, for $k$ big enough, 
\begin{displaymath}
|G_k|\leq \frac{4\sqrt{2}e^{c^2}}{\sqrt{\eta}},
\end{displaymath}
because for $k$ sufficiently big, since $\tilde{q}_k$ converges uniformly to $p$ (after passing to a subsequence), 
\begin{displaymath}
|\tilde{q}_k|\geq  \frac{\sqrt{\eta}}{4\sqrt{2}}.
\end{displaymath}
Furthermore, we get a bound on the derivative of $G_k$
\begin{displaymath}
|G'_k|\leq \frac{64\sqrt{2}e^{3c^2}}{\eta},
\end{displaymath}
so that, we may conclude that $G_k$ has a convergent subsequence (uniformly on compact subsets of $D_c$). By bounding the second derivative of $G_k$, we can even assume that $G'_k$ is convergent in that same subsequence (uniformly on compact subsets of $D_c$). Consider now $\tilde{f}_k$, which is simply $f_k$ rescaled. These functions define a sequence which is uniformly bounded (in fact $f_k$ varies between $\log(\eta)$ and $0$), as is the sequence of derivatives of all orders. In fact
\begin{displaymath}
\tilde{f}_k=\log(|\tilde{p}_k|^2+|\tilde{q}_k|^2)-\frac{k}{2}\log \left( 1+\frac{z\bar{z}}{k}\right),
\end{displaymath}
and, for example,
\begin{displaymath}
\frac{\partial \tilde{f}_k}{\partial z}=\frac{{\tilde{p}_k}'\bar{\tilde{p}}_k+{\tilde{q}_k}'\bar{\tilde{q}}_k}{(|\tilde{p}_k|^2+|\tilde{q}_k|^2)}-\frac{\bar{z}}{2} \left( 1+\frac{z\bar{z}}{k}\right)
\end{displaymath}
which is clearly bounded in $D_c$. This proves that $\tilde{f}_k$ has a subsequence which converges, together with its derivatives of higher order, on compact subsets of $D_c$. Let's denote by $f$ the limit of this subsequence. The function $f$ is smooth and the derivatives of $\tilde{f}_k$ converge to the derivatives of $f$. Next, we show that $G$ and $f$ satisfy a differential relation, that comes from taking the limit of the rescaled version of the differential relation (\ref{PDE}). For this purpose, consider the metric  
\begin{displaymath}
\tilde{\omega}=\frac{dzd\bar{z}}{\left( 1+\frac{|z|^2}{k}\right)^2}
\end{displaymath}
on $D_c$, which is $k$ times the pull back of the Fubini-Study metric on $S^2$ by the rescaling map $\delta_k$ (this map is defined in the z coordinate by $\delta_k(z)=z/\sqrt{k}$). This metric has curvature $1/k$ and tends to a flat metric. Now $G_k=F_k \circ \delta_k$ so that
\begin{displaymath}
\frac{G_k^*\omega_{FS}}{\tilde{\omega}}=\frac{\delta_k^*F_k^*\omega_{FS}}{k\delta^*_k \omega_{FS}}=\frac{1}{k}\frac{F_k^*\omega_{FS}}{\omega_{FS}}\circ \delta_k.
\end{displaymath}
As for the Laplacian of $f_k$ with respect to the Fubini-Study metric, it is given by
\begin{displaymath}
\triangle f_k=(1+|z|^2)^2 (f_k)_{z\bar{z}},
\end{displaymath}
whereas the Laplacian of $\tilde{f}_k$ in the metric $\tilde{\omega}$ is given by
\begin{displaymath}
\triangle \tilde{f}_k=(1+\frac{|z|^2}{k})^2 (\tilde{f}_k)_{z\bar{z}}.
\end{displaymath}
We see that 
\begin{displaymath}
\frac{1}{k}(\triangle f_k)\circ \delta_k=\frac{1}{k}(1+\frac{|z|^2}{k})^2 (f_k)_{z\bar{z}}\circ\delta_k=\triangle \tilde{f}_k.
\end{displaymath}
Then the differential relation (\ref{PDE}) becomes 
\begin{displaymath}
\frac{G_k^*\omega_{FS}}{\tilde{\omega}}=1+\triangle \tilde{f}_k.
\end{displaymath}
Taking the limit of this equation as $k$ tends to infinity (we actually need to consider a subsequence) we get
\begin{equation}\label{PDE_sem_k}
G^*\omega_{FS}=(1+\triangle f)\omega_{\text{flat}},
\end{equation}
where $\omega_{\text{flat}}$ is the flat metric in the disc and the Laplacian is the flat Laplacian. Intuitively, this equation says that, at least when $\eta$ is very close to $1$, the rescaled $F_k$ is trying to be an isometry between flat $\mathbb{C}$ and the sphere. Such isometries do not exist, even locally, so it cannot be that.
We now prove the following lemma:
\begin{lemma}\label{delta_lambda}
Let $c$, $\delta$ and $M$ be positive constants and suppose we have a fixed metric on $D_c$. Then, there is a constant $\lambda=\lambda(\delta,M,c)$, such that, for every smooth function $f$ from the disc $D_c$ to $\mathbb{C}$ satisfying
\begin{itemize}
\item $|f|\leq \delta$,
\item $|d\triangle f|\leq M$,
\end{itemize}
we have $|\triangle f|\leq \lambda$ in $D_{c/4}$.
\end{lemma}
{\bf proof}: Let $z_0$ be a point in $D_{c/4}$ where $\triangle f$ is greater than $\lambda$. We want to see that $\lambda$ big leads to a contradiction. (The same reasoning would lead to a contradiction if we assumed that $\triangle f \leq -\lambda$). Given $\xi$ in $D_{c/4}$, for some $\chi$ in the segment from $z_0$ to $\xi$, we have
\begin{displaymath}
\triangle f (z_0+\xi)- \triangle f (z_0)=d\triangle f (\chi) \xi.
\end{displaymath}
Choosing $\xi$ with norm smaller than $\lambda/2M$, we conclude that $\triangle f (z_0+\xi)$ is greater than $\lambda/2$. Set $\mu$ to be the minimum of $\lambda/2M$ and $c/4$, we have seen
\begin{displaymath}
\triangle f (z_0)\geq \lambda \implies \triangle f (z) \geq \lambda/2, \,\,\,\, z\in B(z_0,\mu).
\end{displaymath}
 Let $r$ denote the distance function to $z_0$ and $\beta$ a smooth function on $\mathbb{R}^+$ that is zero in $]\mu, +\infty[$ and $1$ in $[0,\mu/2]$. Set $g(z)=\beta(r)$, $g$ has compact support contained in $D_c$ and, for any positive number $\epsilon$, $\beta$ can be chosen so as to have $g$ satisfy $\triangle g\leq (8+\epsilon)/\mu^2$. Since $g$ vanishes in a neighborhood of the boundary of $D_c$,
 \begin{displaymath}
 \int_{D_c}g \triangle f = \int_{D_c}f \triangle g. 
 \end{displaymath}
 Now
 \begin{itemize}
 \item the right hand side is equal in norm to
 \begin{displaymath}
 \left | \int_{B(z_0,\mu)}f \triangle g \right| \leq \frac{8+\epsilon}{\mu^2} \delta\pi \mu^2=(8+\epsilon)\pi\delta,
 \end{displaymath}
 \item whereas the left hand side is equal to
 \begin{displaymath}
  \int_{B(z_0,\mu)}g \triangle f \geq  \int_{B(z_0,\mu/2)}g \triangle f \geq\frac{\pi \lambda \mu^2}{8}.
 \end{displaymath}
 \end{itemize}
Therefore we get a contradiction if $\lambda\mu^2> (64+\epsilon)\delta$ (a different $\epsilon$...). Taking
 \begin{displaymath}
 \lambda=\max (7M^{2/3}\delta^{1/3}, \frac{1025\delta}{c^2})
 \end{displaymath}
 we get $|\triangle f (z_0)|\leq \lambda$ for every $z_0$ in $D_{c/2}$.
 
Let us again consider the differential relation (\ref{PDE_sem_k}). The left hand side represents a metric conformal to the flat metric on the disc, with conformal factor $1+\triangle f$, its curvature is therefore
 \begin{displaymath}
- \frac{1}{2} \frac{\triangle \log (1+\triangle f)}{1+\triangle f}.
\end{displaymath}
As for the right hand side, it represents a metric (with singularities maybe), whose curvature is $1$, because it is the pullback of the Fubini-Study metric. We have therefore established the following equality
\begin{equation}\label{equacao_curvatura}
1=- \frac{1}{2} \frac{\triangle \log (1+\triangle f)}{1+\triangle f}.
\end{equation}
We know that $|f_k|\leq \log(1/\eta)$, so that $f$ satisfies the same inequality. From Lemma \ref{delta_lambda}, this implies that $\triangle f$ is small (because $d\triangle f$ is bounded from above by some constant depending only on $\eta$ as we will see ahead) and, so, $ \log (1+\triangle f)$ is small as well. We will show that the derivative of the Laplacian of this quantity is bounded (by a quantity depending only on $\eta$), so, we can apply Lemma \ref{delta_lambda} again, to conclude that $\triangle \log (1+\triangle f)$ is small. Since the denominator of the right hand side of equation (\ref{equacao_curvatura}) is bounded from below ($\triangle f$ is small), then, the right handside itself is small. But it needs to equal $1$, so it can't be small and we arrive at a contradiction. We still need to get an upper bound on the derivatives of $\triangle f$ and $ \log (1+\triangle f)$. From equation (\ref{PDE_sem_k}) we write
\begin{displaymath}
1+\triangle f =\frac{|p'q-pq'|^2}{(|p|^2+|q|^2)^2}.
\end{displaymath}
To estimate $d\triangle f$, we estimate 
\begin{displaymath}
d\frac{|p'q-pq'|^2}{(|p|^2+|q|^2)^2},
\end{displaymath}
on $D_c$, when $c$ is sufficiently small by using the bounds for the norms of $p$ and $q$ (and their derivatives) which, as we remarked before, depend only on $c$. We also use the fact that $|p|^2+|q|^2 \geq \eta$. 
Note that, we can only apply Lemma \ref{delta_lambda} to a $c$ sufficiently small, since it is only on balls with such radii that there is convergence and therefore $f$ is defined. We also need an upper bound for $|d\triangle \log(1+\triangle f)|$. Using the fact that
\begin{displaymath}
\triangle \log(1+\triangle f)=-2\triangle \log(|p|^2+|q|^2),
\end{displaymath}
away from the zeroes of $p'q-qp'$, as well as the inequalities mentioned above for $|(p,q)|$, $|(p',q')|$ and $|(p'',q'')|$, we can obtain such an upper bound depending only on $c$ over $D_c$. Following trough with these calculations and choosing a small enough $c$ in terms of $\eta$, we conclude that 
\begin{displaymath}
\log(1/\eta) \geq \frac{1}{(74708)^9 e^{28+\frac{15}{16}}},
\end{displaymath}
which is of the order $10^{-36}$, giving an upper bound for $\eta$ of the order $e^{-10^{-36}}< 1$. 

\subsection{Uniform distribution and estimated transversality on $S^2$}
Donaldson's construction gives rise to sections which themselves give rise to several submanifolds (the fibers of Lefschetz pencils away from their singularities for example). An important feature of those is that they are asymptotically uniformly distributed.
\begin{theorem}[Donaldson,\cite{do1}]\label{uniformemente_distribuido}
 Let $X$ be a symplectic manifold with an almost-complex structure and $L_k\rightarrow X$ be complex Hermitian line bundles with curvature $i\omega_k/2$, such that we can cover $X$ by balls $B_i=B(p,rk^{-1/2})$ where $L_k$ have local trivializations, $\sigma_i$, such that
\begin{itemize}
\item $0<\alpha\leq |\sigma_i|\leq 1$, for some $\alpha$ independent of $k$,
\item $|\nabla \sigma_i|\leq \beta \sqrt{k}$.
\end{itemize}
Suppose further, that $\{s_k\}$ is a sequence of sections of $L_k$ for which there is a sequence ${a_k}$ of numbers, bounded away from zero, such that  
\begin{itemize}
\item $|s_k|\leq aa_k$, $|\nabla s_k|\leq b \sqrt{k}a_k$,
\item $|\bar \partial s_k|=O(a_k)$ ,
\item $|s_k|<a_k \implies |\nabla s_k|\geq 2\frac{\beta}{\alpha}a_k\sqrt{k}$.
\end{itemize}
Then, the sets $s_k^{-1}(0)$ are asymptotically uniformly distributed, with the following estimate:
\begin{displaymath}
\mid \int_{s_k^{-1}(0)}\psi -\frac{1}{2\pi}\int_X \omega_k\wedge \psi \mid \leq \sqrt{k}|d\psi|_{\infty},
\end{displaymath}
for all $\psi$, $2n-2$ form on $X$.
\end{theorem}
The statement here is slight modification of the corresponding statement in \cite{do1}. The proof is almost the same, so we won't reproduce it here.
We can apply this theorem to the bundles $L^k$ over $X$, and the sequence of sections $s_0$ coming from Theorem \ref{exist_lapiz}. The sequence ${a_k}$ can be chosen to be constant equal to $\eta$. In this way we recover the result of \cite{do1}. We can also apply the theorem to the sections $s_0-bs_1$, when $|b|$ is small (with respect to $\eta$) and $s_0$ satisfies the conditions in Theorem \ref{exist_lapiz}. Another interesting result of the same type can be obtained by applying Theorem \ref{uniformemente_distribuido} to the bundles $L^k\otimes L^k \otimes T^*X$, and their sections $s_1\partial s_0-s_0\partial s_1$, when $s_1$ and $s_2$ come from Theorem \ref{exist_lapiz} and $\eta$ is small. Here we take $a_k=\epsilon \sqrt{k}$, for some $\epsilon$ to be made explicit shortly. Let's first check that the sections satisfy the hypotheses in the theorem. We know
\begin{displaymath}
\nabla \left(\frac{s_1}{s_0}\right)=\frac{s_1\nabla s_0-s_0\nabla s_1}{s_0^2}
\end{displaymath}
and
\begin{displaymath}
\nabla(s_1\partial s_0-s_0\partial s_1)=s_0^2\nabla \nabla \left(\frac{s_1}{s_0}\right)-(s_1\nabla s_0-s_0\nabla s_1)\frac{2\nabla s_0}{s_0}.
\end{displaymath}
Now, 
\begin{displaymath}
|s_1\nabla s_0-s_0\nabla s_1|\leq \frac{\eta^2\sqrt{k}}{2\sqrt{2}} \implies |s_0|\geq \frac{\eta^2}{2\sqrt{2}C}
\end{displaymath}
(where $C$ comes from $|\partial s_1|\leq C\sqrt{k}$). This is because, 
\begin{displaymath}
|s_0|\leq \frac{\eta^2}{2\sqrt{2} C} \implies |s_0|\leq \eta
\end{displaymath}
which in turn implies that $|\nabla s_0|\geq \eta\sqrt{k}$. On the other hand, because $ {\eta^2}/{2\sqrt{2}C}$ is $O(\eta/\sqrt{2})$, this also implies that $|s_1|\geq \eta/\sqrt{2}$, so that, when $|s_0|\leq \eta^2/(2\sqrt{2}C)$, 
\begin{displaymath}
|s_1\nabla s_0-s_0\nabla s_1|\geq \frac{\sqrt{k}\eta^2}{\sqrt{2}}-\frac{\eta^2}{2\sqrt{2}C}C\sqrt{k}=\frac{\sqrt{k}\eta^2}{2\sqrt{2}}.
\end{displaymath}
Assume that $|s_1\nabla s_0-s_0\nabla s_1|\leq \eta^8\sqrt{k}$, then (if $\eta$ is small enough to make $\eta^8 \leq  \frac{\eta^2}{2\sqrt{2}}$) we have, 
\begin{displaymath}
\left|\nabla \left(\frac{s_1}{s_0}\right)\right|=\frac{|s_1\nabla s_0-s_0\nabla s_1|}{|s_0^2|}\leq a\sqrt{k}\eta^6
\end{displaymath} 
therefore, $|\partial \partial ({s_1}/{s_0})|\geq k\eta$. Then, 
\begin{displaymath}
|\nabla(s_1\nabla s_0-s_0\nabla s_1)|\geq a\eta^4k\eta-\frac{\sqrt{k}\eta^82C\sqrt{k}}{\eta^2}\geq \eta^8
\end{displaymath}
when $\eta$ is sufficiently small. We have shown that the sequence of sections of $L^k\otimes L^k \otimes T^*X$, $s_1\nabla s_0-s_0\nabla s_1$ is $\epsilon=\eta^8$ transverse to zero. The zeroes of this sequence of sections form the set of branch points of the Lefschetz pencil. We then conclude:
\begin{prop}
The set of branch points of a Lefschetz pencil coming from Theorem \ref{exist_lapiz_holo} is asymptotically uniformly distributed. 
\end{prop}
On $S^2$, the uniform distribution assumes a particularly simple expression since the fibers of the pencil are sets of $k$ points, which we denote by $\{x_i^k\}$, satisfying
\begin{equation}\label{estimativa_uniforme}
\left| \frac{1}{k}\sum_{i=1}^{k} f(x_i^k)-\frac{1}{|S^2|}\int_{S^2} f\right| \leq C\frac{|df|_{\infty}}{\sqrt{k}},
\end{equation}
for all $\mathcal{C}^{\infty}$ functions $f$.
The branch point set is simply a set of $2k-2$ points, $\{p_i^k\}$. It satisfies
\begin{displaymath}
\left| \frac{1}{2k-2}\sum_{i=1}^{2k-2} f(p_i^k)-\frac{1}{|S^2|}\int_{S^2} f\right| \leq C\frac{|df|_{\infty}}{\sqrt{k}},
\end{displaymath}
for all $\mathcal{C}^{\infty}$ functions $f$.

Asymptotically uniformly distributed sets of points on $S^2$ are important in many areas of mathematics and are the object of a lot of work in potential theory (see \cite{rsz} and analytic number theory (see \cite{lps}, \cite{s}). Therefore, these properties are very helpful in trying to explicitly build the sequence of polynomials $(p_k,q_k)$. On the other hand, they are certainly not enough to characterize these polynomials. In fact, asymptotic uniform distributiveness with estimate (\ref{estimativa_uniforme}) would allow for a point with multiplicity, whereas Donaldson's construction clearly does not. 

A natural question to ask at this point is: How about the fibers above points close to $\infty$ on $S^2$? We have not shown that these are asymptotically uniformly distributed since our method failed for them. Next, we show that in fact they are, in the case of $S^2$ (although in a slightly different sense). What is more surprising, we show (using the methods developed in the previous section) that the condition of transversality of the pair is enough to ensure the asymptotic uniform distribution of all the fibers and also of the branch points, for pencils in $S^2$. As before, we identify the sections of $\mathcal{O}(k)$ with homogeneous polynomials and write $||.||$ for their norm.

{\bf proof of Theorem \ref{ud}}: We will start by showing that the zero set of $p_k$ is uniformly distributed in the sense that, if $x^k_i$ denote the points in this set for $i=1, \cdots k$, counted with multiplicity, and $f$ is a $\mathcal{C}^2$ function on $S^2$, then
\begin{equation}
\left| \frac{1}{k}\sum_{i=1}^k f(x^k_i)-\frac{1}{|S^2|}\int_{S^2} f \right|\leq \frac{C||\triangle f||_{\infty}}{k},
\end{equation}
\begin{enumerate}
\item The first step is to show that 
\begin{displaymath}
\int_{S^2} \log ||p_k||^2
\end{displaymath}
is bounded independently of $k$. We will do this by contradiction. Let's assume that it is not. For each $k$, it is easy to find a partition of $S^2$ by $k$ sets $A_i$ such that $\text{diam}(A_i)=c/\sqrt{k}$ ($c$ is independent of $k$). By assumption, there is a subsequence $n_k$, such that $\int_{S^2} \log ||p_{n_k}||^2$ tends to $-\infty$ as $k$ tends to $+\infty$ (note that $\log ||p_{n_k}||^2\leq 0$) . Now take $i=i(k)$ to be the index of the element in the partition of $S^2$ were the integral of $\log ||p_{n_k}||^2$ is the smallest. We have
\begin{displaymath}
n_k\int_{A_i}  \log ||p_{n_k}||^2\leq \int_{S^2} \log ||p_{n_k}||^2,
\end{displaymath}
so that the left hand side tends to $-\infty$ as well. But $A_i\subset B(z_k,c/\sqrt{k})$ so
\begin{displaymath}
n_k\int_{B(z_k,c/\sqrt{n_k})}  \log ||p_{n_k}||^2
\end{displaymath}
tends to $-\infty$ as well. Suppose we knew that $z_k$ was zero for all $k$. Then, by using the rescaled coordinate chart around zero which we described in the proof of Proposition \ref{upper_bound}, we would have for this integral
\begin{displaymath}
\int_{B(0,c)}\log \left( \frac{|\tilde{p}_{n_k}|^2}{\left(1+\frac{|z|^2}{n_k}\right)^{n_k}}\right) \frac{dzd\bar{z}}{\left(1+\frac{|z|^2}{n_k}\right)^{2}}
\end{displaymath}
(were $\tilde{p}_k$ is simply $p_k$ in the rescaled coordinates as in the proof of Proposition \ref{upper_bound}). We showed before that $\{\tilde{p}_{n_k}\}$ is bounded and has bounded derivatives in $B(0,c)$ and therefore has a convergent subsequence. The limit, $p$, is a holomorphic function on $B(0,c)$. We will check ahead that it cannot be identically zero and therefore has a finite number of zeroes, with finite multiplicity, so that $\log |p|^2$ is integrable. The integral above converges for that subsequence to the finite number
\begin{displaymath}
\int_{B(0,c)}\log \left( \frac{|{p}|^2}{e^{|z|^2}} \right) {dzd\bar{z}}
\end{displaymath}
and we get a contradiction. Let's try to reproduce this argument for a general sequence of points $z_k$. Given any point $p$ on $S^2$ there is a map $\phi:S^2 \rightarrow S^2$ such that 
\begin{enumerate}
\item $\phi^* \omega_{FS}= \omega_{FS}$,
\item $\phi$ takes the north pole to $p$,
\item the inverse of $\phi$, $\psi$ is such that $|\psi^* \mathbf{w}|=|\mathbf{w}|\circ \psi$ were $[\mathbf{z}:\mathbf{w}]$ are homogeneous coordinates on $S^2$, $\mathbf{w}$ denotes the same section of $\mathcal{O}(1)$ as before and $|.|$ is the usual norm in $\mathcal{O}(1)$. 
\end{enumerate}
Such a map can be taken to simply be a rotation. Now consider the composition of this map for $p=z_k$ with the inverse of the usual chart centered at the north pole, $l$, $\phi \circ l:\mathbb{C}\rightarrow S^2$ (we omit the dependence on $k$). Its inverse defines a chart, centered at the point $z_k$. Write $p_k$ as a holomorphic function on $\mathbb{C}$ times the section of $\mathcal{O}(k)$, $\psi^* \mathbf{w}^k$. Let $f_k$ denote this holomorphic function composed with $\phi \circ l$. We have, 
\begin{displaymath}
|\psi^*\mathbf{w}^k|\circ \phi \circ l=|\mathbf{w}^k|\circ \psi \circ \phi \circ l=|\mathbf{w}^k|\circ l =\frac{1}{\left(1+{|z|^2}\right)^{\frac{k}{2}}}.
\end{displaymath}
As a consequence
\begin{displaymath}
|f_k|^2\leq {\left(1+{|z|^2}\right)^{k}}.
\end{displaymath}
The same can be done for $q_k$, so that we obtain two sequences of functions $f_k$ and $g_k$ satisfying 
\begin{displaymath}
\eta \left(1+{|z|^2}\right)^{k} \leq |f_k|^2+|g_k|^2\leq {\left(1+{|z|^2}\right)^{k}}.
\end{displaymath}
If we rescale the $z$ coordinate in $\mathbb{C}$ by $1/\sqrt{k}$ as we did in the proof of Proposition \ref{upper_bound}, then, we get functions $\tilde{f}_k$ and $\tilde{g}_k$ which are bounded on $B(0,c)$, independently of $k$, for any fixed $c$, and whose derivatives are also bounded there. They then have convergent subsequences (in the $\mathcal{C}^\infty$ norm over $B(0,c)$). In fact, we will use that $\tilde{f}_{n_k}$ and $\tilde{g}_{n_k}$ have convergent subsequences. We call $p$ and $q$ their limits. As before,
\begin{displaymath}
n_k\int_{B(z_k,c/\sqrt{n_k})}  \log ||p_{n_k}||^2=\int_{B(0,c)}\log \left( \frac{|\tilde{f}_{n_k}|^2}{\left(1+\frac{|z|^2}{n_k}\right)^{n_k}}\right) \frac{dzd\bar{z}}{\left(1+\frac{|z|^2}{n_k}\right)^{2}}.
\end{displaymath}
The sequences $\tilde{f}_{n_k}$ and $\tilde{g}_{n_k}$ are subconvergent on compact subsets of $\mathbb{C}$ to $p$ and $q$. These satisfy
\begin{displaymath}
\eta e^{|z|^2}\leq |p|^2+|q|^2\leq e^{|z|^2}.
\end{displaymath}
This inequality implies that $p$ cannot be identically zero. If it were, there would be a holomorphic function $q$ satisfying
\begin{displaymath}
\eta e^{|z|^2}\leq |q|^2\leq e^{|z|^2}
\end{displaymath}
so that $q$ could not vanish in $\mathbb{C}$. Twice the real part of its log, $u$, would satisfy
\begin{displaymath}
\log \eta+|z|^2\leq u\leq |z|^2.
\end{displaymath}
The function $u-|z|^2$ would be bounded on $\mathbb{C}$ and therefore would have a minimum, but its Laplacian is $-2\leq0$, so it could not have a minimum. For $k$ sufficiently big, on $B(0,c)$, since (a subsequence of) $\{\tilde{f}_{n_k}\}$ converges uniformly to $p$ 
\begin{displaymath}
\log \frac{|\tilde{f}_{n_k}|^2}{\left(1+\frac{|z|^2}{n_k}\right)^{n_k}}\geq \frac{1}{2}\log \frac{|p|^2}{e^{|z|^2}},
\end{displaymath}
and therefore
\begin{displaymath}
0\geq n_k\int_{B(z_k,c/\sqrt{n_k})} log ||p_{n_k}||^2 \geq \frac{1}{2}\int_{B(0,c)} \log \frac{|p|^2}{e^{|z|^2}},
\end{displaymath}
for that same subsequence. This is a contradiction since we are assuming the left hand side integral converges to $-\infty$.
\item Given the zeroes of $p_k$ as a function on $S^2$, $x_i$ (we omit the $k$ dependence here), for each of these, consider the function on $S^2$ given by 
\begin{displaymath}
\log \left( d^2(x,x_i)\frac{e}{4}\right),
\end{displaymath}
where $d$ is the chordal distance in $S^2$, i.e., the usual distance in $\mathbb{R}^3$. The Laplacian of this function is simply $4\pi\delta_{x_i}-1$ and its integral is zero (hence the $e/4$). Let $z_i$ be the complex coordinates of the points $x_i$, trough the usual chart centered at the north pole. Now consider the homogeneous polynomial 
\begin{displaymath}
\hat{p_k}=\prod_{i}\left( \frac{e^{1/2}(z-z_i)}{\left( 1+|z|^2\right)^{1/2}} \right) \mathbf{w}^k.
\end{displaymath}
As a section of $\mathcal{O}(1)$, the norm square of this is 
\begin{displaymath}
\prod_{i} \left( d^2(x,x_i)\frac{e}{4}\right).
\end{displaymath}
This is because, given two points on the sphere, $x$ and $y$ with complex coordinates $z$ and $w$
\begin{displaymath}
d^2(x,y)=\frac{4|z-w|^2}{(1+|z|^2)(1+|w|^2)}.
\end{displaymath}
Such a $\hat{p_k}$ satisfies $\int \log||\hat{p}_k||^2=0$. The polynomials $p_k$ and $\hat{p_k}$ are related by $\hat{p_k}=c_k p_k$, for some sequence of numbers $c_k$. The fact that $\int \log||{p}_k||^2$ is bounded independently of $k$, implies that this sequence is bounded. 
\item Finally, let $f$ be a $\mathcal{C}^2$ function on $S^2$ and consider the difference
\begin{displaymath}
\frac{\sum f(x_i)}{k}-\frac{1}{4\pi}\int f.
\end{displaymath}
This can alternatively be written as 
\begin{displaymath}
\frac{1}{4\pi}\int f\frac{\sum \left( 4\pi \delta_{x_i}-1\right)}{k},
\end{displaymath}
which, in view of a previous remark, is
\begin{displaymath}
\int f\frac{\sum \triangle \log \left( d^2(x,x_i)\frac{e}{4}\right)}{k}=\frac{1}{k}\int f\triangle \log||\hat{p}_k||^2.
\end{displaymath}
Integrating by parts, we get,
\begin{displaymath}
\frac{1}{k}\int \triangle f \log||\hat{p}_k||^2.
\end{displaymath}
To prove the proposition, it is then sufficient to show that $\int \left | \log||\hat{p}_k||^2\right |$ is bounded, independently of $k$. This integral is simply bounded by
\begin{displaymath}
4\pi\log |c_k|^2-\int \log ||p_k||^2=-2\int \log ||p_k||^2,
\end{displaymath}
which is bounded, as we have already seen.
\end{enumerate}
In the same way, to prove that any fiber is uniformly distributed in the sense described previously, it is enough to see that the integral
\begin{displaymath}
\int_{S^2} \left| \log ||p_k+\lambda q_k||^2\right|
\end{displaymath}
is bounded independently of $k$ for any $\lambda$. Assume that it is not. By the same reasoning as before, we can conclude that, for a sequence of complex numbers $z_k$,
\begin{displaymath}
\int_{B(0,c)} \left| \log \frac{| \tilde{f}_{n_k} +\lambda \tilde{g}_{n_k} |^2}{\left( 1+\frac{|z|^2}{n_k} \right)^{n_k}} \right| \frac{dzd\bar{z}}{\left( 1+\frac{|z|^2}{n_k} \right)^{2}}\rightarrow \infty,
\end{displaymath}
where $\tilde{f}_k$ and $\tilde{g}_k$ are defined as before using the $z_k$. But, again as before, $\{\tilde{f}_{n_k}\}$ and $\{\tilde{g}_{n_k}\}$ subconverge uniformly on compact subsets of $\mathbb{C}$, so that, for each $\lambda$, $\{\tilde{f}_{n_k}+\lambda \tilde{g}_{n_k\}}$ subconverges to $p+\lambda q$, uniformly on compact subsets of $\mathbb{C}$. To be able to apply the same reasoning as before, we need only check that $p+\lambda q$ is not identically zero. If it were, there would be a holomorphic function, $g$ satisfying
\begin{displaymath}
\frac{\eta}{(1+|\lambda|^2)}e^{|z|^2}\leq |g|^2 \leq \frac{1}{(1+|\lambda|^2)}e^{|z|^2}, 
\end{displaymath}
which is impossible. As for the branch points, we need to show that 
\begin{equation}\label{integral_branch_points}
\int_{S^2} \left| \log \frac{||p_k\nabla q_k -q_k\nabla p_k||^2}{k} \right|
\end{equation}
is bounded independently of $k$. Note that $p_k\nabla q_k -q_k\nabla p_k$ is a section of $T^*S^2\otimes \mathcal{O}(k)\otimes \mathcal{O}(k)$, which is simply $\mathcal{O}(2k-2)$ and $||.||$ refers to the usual norm in that bundle. Letting $f_k$ and $g_k$ be the representations of $p_k$ and $q_k$ respectively in the trivialization $\psi^*\mathbf{w}^k$ used before, we can write
\begin{displaymath}
 p_k\nabla q_k -q_k\nabla p_k=(f_kg'_k-f'_kg_k)dz\otimes \psi^*\mathbf{w}^k \otimes \psi^*\mathbf{w}^k.
 \end{displaymath}
 The norm square of this section is given by
 \begin{displaymath}
\frac{ |f_kg'_k-f'_kg_k|^2}{\left( 1+|z|^2\right)^{2k}}. \frac{1}{2}\left( 1+|z|^2\right)^2=\frac{ |f_kg'_k-f'_kg_k|^2}{2\left( 1+|z|^2\right)^{2k-2}}.
\end{displaymath}
By rescaling the coordinate $z$, we get
\begin{displaymath}
\tilde{f'_k}(z)=\frac{1}{\sqrt{k}}f'_k\left(\frac{z}{\sqrt{z}}\right).
\end{displaymath}
Assuming that the sequence of integrals in (\ref{integral_branch_points}) is not bounded, we conclude that there must be complex numbers $z_k$ such that
\begin{displaymath}
\int_{B(0,c)} \left| \log \frac{|\tilde{f}_{n_k}\tilde{g}'_{n_k}-\tilde{f}'_{n_k}\tilde{g}_{n_k}|^2}
{2\left(1+\frac{|z|^2}{n_k}\right)^{2n_k-2}}\right| \frac{dzd\bar{z}}{\left(1+\frac{|z|^2}{n_k}\right)^{2}} 
\end{displaymath}
tends to infinity. The sequence $\{\tilde{f}_{n_k}\tilde{g}'_{n_k}-\tilde{f}'_{n_k}\tilde{g}_{n_k}\}$ is uniformly subconvergent on compact subsets of $\mathbb{C}$ to $pq'-p'q$. We will be done if we check that this function cannot be identically zero. If it was, then $p/q$ would be constant, which is impossible, as we saw above (note that $q$ is not identically zero either). 

\section{Optimally distributed points on $S^2$}

\subsection{Logarithmic equilibrium points on $S^2$}
The problem of finding an optimal way of distributing points on the sphere has been studied in several areas of mathematics. An instance of this problem appears in potential theory. Namely, let $\{x_i\}$ be a set of $k$ points on $S^2$. One can consider the function $\Pi_{i<j} d(x_i,x_j)$, where $d(x,y)$ denotes the distance in $\mathbb{R}^3$, between any pair of points $x$ and $y$ of the sphere. This function achieves its minimum when two of the points of the set coincide. The question is: for which sets of points does it attain its maximum? Maximizing $\Pi_{i<j} d(x_i,x_j)$ is the same as minimizing the logarithmic potential, 
\begin{displaymath}
\sum_{i<j} \log\frac{1}{d(x_i,x_j)}
\end{displaymath}
and the points which achieve the minimum are called logarithmic equilibrium points. It is not known what this (these) configuration(s) is (are), but some things are know (or conjectured) about the minimum value of the logarithmic potential. As a motivation for what follows, let us first see how this problem is related to that of determining the best $\eta$ for which
\begin{displaymath}
 ||p_k||\leq \eta \implies ||\nabla p_k||\geq \eta\sqrt{k}. 
\end{displaymath} 
Given a set of points $\{x_i\}$, by considering the complex coordinates  ${z_i}$ of these points trough stereographic projection, we can form a complex polynomial of degree $k$, whose zeroes are the $z_i$'s. It is determined up to a multiplicative constant. Consider as before 
\begin{displaymath}
p_k(z)=e^{k/2}\prod_i \frac{(z-z_i)}{{(1+|z_i|^2)}^{1/2}}.
\end{displaymath}
Let $x$ denote the point on $S^2$ whose complex coordinate is  $z$. As a section of $\mathcal{O}(k)$ (i.e., as $p_k\mathbf{w}^k$), the norm of $p_k$ at $x$ is
\begin{equation}\label{norma_poly}
|| p_k || ^2=\prod_i \frac{e}{4}d^2(x,x_i).
\end{equation}
There are two things to note here. First recall that, given two points $x$ and $y$ on $S^2$, if you denote their complex coordinates by $z$ and $w$ respectively,
\begin{displaymath}
d(x,y)=\frac{2|z-w|}{(1+|z|^2)^{1/2}(1+|w|^2)^{1/2}}.
\end{displaymath}
Second, using the metric on $S^2$ given by $2idzd\bar{z}/(1+|z|^2)^2$, whose volume is $4\pi$,
\begin{displaymath}
\frac{1}{4\pi}\int_{S^2}\log {d(x,y)}dx=\frac{1}{2}\log\frac{4}{e}
\end{displaymath}
which corresponds to the fact that the integral of the logarithm of the norm square of $p_k$ is zero. 

In order for $\{p_k\}$ to define an $\eta$ transverse sequence of sections of $\mathcal{O}(k)$ we need $|| \nabla p_k(z_i)|| $ to be as big as possible (greater than $\eta\sqrt{k}$ at least), where, again, $p_k$ is thought of as section of  $\mathcal{O}(k)$. We have
\begin{displaymath}
|| \nabla p_k || (z_i)=\frac{1}{\sqrt{2}}\left( \frac{e}{4}\right)^{k/2}\Pi_{i\ne j} d(x_i,x_j).
\end{displaymath}
Here we used the norm of $dz$ which is $(1+|z|^2)/2\sqrt{2}$. Having $\prod_i || \nabla p_k || (z_i)$ big would be to our advantage. But this product is simply a constant (depending only on $k$) times $\Pi_{i<j} d(x_i,x_j)$. The problem of maximizing this quantity is therefore related to our own problem. In fact, we will see ahead, that a good approximation of these logarithmic points gives rise to a pairs of polynomials which satisfy Donaldson's constraints, as in Proposition \ref{principal}, or at least seem to, experimentally. The fact that, in some sense, this distribution is "better" than its approximation leads to the conjecture that these polynomials define sections that do, in fact, satisfy the properties in Theorem \ref{exist_lapiz_holo}. In particular, we conjecture that $||p_k||$ is bounded by a constant independent of $k$. Assuming that  $\liminf \max||p_k||$ is not zero, we can then find an upper bound for $\eta$ in terms of $\liminf \max||p_k||=\alpha$. We will show that, if $\eta$ is such that
\begin{displaymath}
||p_k||\leq \eta\implies ||\nabla p_k||\geq \eta\sqrt{k},
\end{displaymath}
then 
\begin{displaymath}   
\eta \leq \frac{\sqrt{\pi }}{\sqrt {e}}(1-e^{-a})^{b/2}, 
\end{displaymath}
where
\begin{displaymath}
a=\frac{2\sqrt{2\pi}} {\sqrt{27}} (\sqrt{2\pi}+\sqrt{2\pi+\sqrt{27}}) ,\,\,\,\, b=\frac{\sqrt{2\pi+\sqrt{27}}-\sqrt{2\pi}} {\sqrt{2\pi+\sqrt{27}}+\sqrt{2\pi}}.
\end{displaymath}
So, the real upper bound for Donaldson's $\eta$, is this quantity divided by $\liminf \max||p_k||$. (This is really just the $\eta$ corresponding to the first condition). Suppose that the sequence of polynomials $p_k$ defines an $\eta$ transverse sequence. At $z_i$, we must have $||\nabla p_k||(z_i)\geq \eta \sqrt{k}$. Since $p_k$ vanishes at $z_i$, we have
\begin{displaymath}
||\nabla p_k(z_i)||=\frac{1}{\sqrt{2}}\prod_{i\ne j} \left(\frac{e}{4}\right)^{1/2} d(x_i,x_j),
\end{displaymath}
so that 
\begin{displaymath}
\prod_i ||\nabla p_k (z_i)||=\left(\frac{1}{2}\right)^{k/2} \prod_{i< j} \frac{e}{4}d(x_i,x_j)^2,
\end{displaymath}
which is related to the logarithmic potential. There are know estimates for this potential, namely:
\begin{theorem}[Rakhmanov, Saff, Zhou, \cite{rsz}]
The logarithmic equilibria points, ${x_i}$, satisfy the following estimate
\begin{displaymath}
\prod_{i< j} \frac{e}{4}d(x_i,x_j)^2\leq \left(\frac{4\pi}{2e}\right)^{k/2}k^{k/2} (1-e^{-a}+\epsilon_k)^{bk/2},
\end{displaymath}
where $a$ and $b$ are as above and $\epsilon_k \rightarrow 0$ as $k\rightarrow \infty$.
\end{theorem}
On the other hand, if we are assuming $\eta$ transversality, we have 
\begin{displaymath}
\prod_i ||\nabla p_k (z_i)||\geq \eta^k k^{k/2}.
\end{displaymath}
Combining these two estimates we obtain the desired result.

\subsection{Generalized spiral points}\label{pontos_espirais}
Even though the exact configuration of logarithmic points is not know, in \cite{rsz}, Rakhmanov, Saff and Zhou describe a set of points they call generalized spiral that give extremely good numerical estimates for the maximum of $\prod_{i\leq j}d(x_i,x_j)$. We will use a modification of these points to construct a pair of polynomials, for which we can verify Donaldson's conditions experimentally. We also have some partial results towards proving that these indeed satisfy the conditions in Theorem \ref{exist_lapiz_holo} (and in particular Proposition \ref{principal}).
Consider the points with cylindrical coordinates $(h,\theta)$ on $S^2$ given by:
\begin{displaymath}
h_i=-1+\frac{2i-1}{k} \;\;\;\;\; i=1, \cdots k
\end{displaymath}
\begin{displaymath}
\theta_0=0,\;\;\theta_i=\theta_{i-1}+\frac{3.6}{\sqrt{k}\sqrt{1-h_i^2}}
\;mod 2\pi\;\;\;\;\; i=1, \cdots k.
\end{displaymath}
Let $z_i$ be the complex coordinates of these points. As before, we consider the polynomials
\begin{displaymath}
p_k=e^{k/2}\prod_i \frac{(z-z_i)}{{(1+|z_i|^2)}^{1/2}},
\end{displaymath}
and \begin{displaymath}
q_k=e^{k/2}\prod_i \frac{(z+z_i)}{{(1+|z_i|^2)}^{1/2}}.
\end{displaymath}
\begin{conj}\label{conjectura_principal}
Up to dividing by a constant that will make 
\begin{displaymath}
\max \left( || p_k ||^2+ || q_k ||^2 \right)
\end{displaymath}
 equal to $1$, the sections $(p_k\mathbf{w}^k,q_k\mathbf{w}^k)$ satisfy the conditions in Theorem \ref{exist_lapiz_holo}. 
\end{conj}
We have verified this experimentally by verifying that:
\begin{itemize}
\item $\max \left( || p_k ||^2+ || q_k ||^2 \right)$ is indeed bounded independently of $k$.
The experimental values found for this maximum are
\begin{center}
\begin{tabular}{|c|c|c|c|c|c|c|c|} \hline
  $k$&50 &100 &150 &170 &180 &190 &
       200\\
\hline
  $\max$&4.7828 &4.8272 &4.8368 &4.8432 &4.8373 &4.8364 &4.8460\\
\hline
\end{tabular}
\end{center}
\item $\min \left( || p_k ||^2+ || q_k ||^2 \right)$ in bounded from below independently of $k$. The experimental values found for the minimum are
\begin{center}
\begin{tabular}{|c|c|c|c|c|c|c|c|} \hline
  $k$&50 &100 &150 &170 &180 &190 &
       200\\
\hline
  $min$&0.7194 &0.7093 &0.7085 &0.7085 &0.7125 &0.7089 &0.7085\\
\hline
\end{tabular}
\end{center}
\item $\min || \nabla p_k ||(z_i)$ is bounded independently of $k$ and the same for $q_k$. Below is a table of the values of the minimum of $||\nabla p_k||^2(z)$ (normalized by $2/\sqrt{k}$) taken among the $z_i$'s
\begin{center}
\begin{tabular}{|c|c|c|c|c|c|c|c|} \hline
  $k$&100 &200 &500 &700 &900 &1000\\
\hline
  $\min\nabla$&1.6963 &1.6998 &1.7020 &1.7024 &1.7026 &1.7027\\
\hline
\end{tabular}
\end{center}
\end{itemize}
These values give an experimental $\eta$ close to $0.15$. We can also show that the sequences  $|| p_k ({z}/{\sqrt{k}})||^2$ for $|z|\leq 1$ are bounded (the same holds true for the corresponding $q$ sequence). We start by showing that $|| p_k (0)||^2$ is bounded. Now $d(x,x_i)^2=2-2\langle x,x_i\rangle$, because $x_i$ and $x$ have norm 1. If $x$ has coordinates $(0,0,-1)$, then $\langle x,x_i\rangle$ is simply minus the last coordinate of $x_i$, that is $-h_i$. Equation (\ref{norma_poly}) becomes
\begin{displaymath}
|| p_k(0) ||^2=\prod \frac{e}{2}(1+h_i)=\prod_{i=1}^k \frac{e(2i-1)}{2k}=2\frac{e^k(2k-1)!}{4^kk^k(k-1)!}.
\end{displaymath}
We use Stirling's approximation for factorial, $k!\sim \sqrt{2\pi k}k^ke^{-k}$, to see that $|| p_k(0) ||^2$ actually converges to $\sqrt{2}$. As for $|| p_k(z/\sqrt{k}) ||^2$,
\begin{displaymath}
|| p_k(\frac{z}{\sqrt{k}}) ||^2=e^{k}\prod \frac{|z_i|^2}{1+|z_i|^2}\prod \left| \frac{z}{z_i\sqrt{k}}-1 \right|.
\end{displaymath}
The first product is simply $|| p_k(0) ||^2$, which we know is bounded. As for the second product, if one considers the inequality $|\log (1-z)|\leq 4|z|$, which holds true for all $z$ with $|z|\leq 1$, then one can see that,
\begin{displaymath}
\left| \log\left| \prod \left( 1-\frac{z}{z_i\sqrt{k}} \right)\right|\right|\leq 4\sum \frac{|z|}{\sqrt{k}|z_i|}.
\end{displaymath}
In cylindrical coordinates, 
\begin{displaymath}
|z_i|=\sqrt{\frac{1-h_i}{1+h_i}}=\sqrt{\frac{2k-2i-1}{2i-1}}.
\end{displaymath}
Now
\begin{displaymath}
\sum_{i=1}^k  \frac{1}{\sqrt{k}|z_i|}=\sum \frac{\sqrt{2i-1}}{\sqrt{k}}\frac{1}{\sqrt{2k-2i+1}}\leq \sqrt{2}\sum_{i=1}^k \frac{1}{\sqrt{i}},
\end{displaymath}
but $\sum 1/\sqrt{i}$ is divergent. If instead of having a single point with given $h_i$, we had 3 (that differed from each other by multiplication by a cube root of unity in complex coordinates), then our initial product would be 
\begin{displaymath}
\prod \left( 1-z^3/(z_i\sqrt{k})^3 \right).
\end{displaymath}
The norm of the logarithm of this product is bounded by a constant times $\sum 1/i^{3/2}$. This new set of points would have $3k$ elements instead of $k$. This would just mean constructing a subsequence of a Donaldson sequence. We have thus showed the following:
\begin{prop}
Consider the set of $3k$ points on $S^2$, whose $h$ coordinates are given by
\begin{displaymath}
h_{3i-2}=h_{3i-1}=h_{3i}=-1+\frac{2i-1}{k} \;\;\;\;\; i=1, \cdots k
\end{displaymath}
and whose $\theta$ coordinates are
\begin{displaymath}
\theta_{3i-2}=\theta_{3i-5}+\frac{3.6}{\sqrt{k}\sqrt{1-h_i^2}}
,\;\;\;\; \theta_{3i-1}=\theta_{3i-2}+\frac{\pi}{3}, \;\;\;\;\; \theta_{3i}=\theta_{3i-2}+\frac{2\pi}{3},\;\;\;\;\; 
\end{displaymath}
defined $\mod 2\pi$ for $i=1, \cdots k$ (where $\theta_0=0$). Then, letting $z_i$ be the complex coordinates of these points trough stereographic projection, and
\begin{displaymath}
p_k=e^{3k/2}\prod_i \frac{(z-z_i)}{{(1+|z_i|^2)}^{1/2}},
\end{displaymath}
$p_k(z/\sqrt{k})$ is bounded, in any fixed neighborhood of 0.
\end{prop}

We cannot do this calculation for the other points of $S^2$, in fact we cannot prove that $||p_k||^2$ is bounded. This is the main ingredient missing. We should remark that, the experimental data shows that the maximum of $||p_k||^2$ occurs close to the north pole, so that this calculations is more relevant to prove what we want, than what it may seem. 

Supposing that the sequence of functions $p_k/(1+|z|^2)^{k/2}$, rescaled, converges, then it is simple to determine what type of limit they should have. To see this, consider one of the spiral points with coordinates $(h_i,\theta_i)$. Suppose we take a small neighborhood, of radius of the order $k^{-1/2}$, around that point. Will there be other spiral points in the neighborhood? The next spiral point, close to $(h_i,\theta_i)$, will appear after a variation in $\theta$ which is approximately $2\pi$, that is, after approximately  
\begin{displaymath}
j=2\pi\sqrt{k}\frac{\sqrt{1-h_i^2}}{3.6}
\end{displaymath}
steps.  Then, the $h$ will be approximately 
\begin{displaymath}
h_i+\frac{4\pi\sqrt{1-h_i^2}}{3.6\sqrt{k}}.
\end{displaymath}
When one considers what this configuration of points looks like, as $k$ tends to $\infty$, one sees that the limit of $p_k/(1+|z|^2)^{k/2}$, rescaled by $\sqrt{k}$ around any point distinct from the poles, can only be a function whose zeroes lie in $\{n+\lambda_m+im, m,n \in \mathbb{Z}\}$, where $0 \leq \lambda_n \leq 1$. We may assume that $\lambda_0=0$. The first thing to check is that such a function actually exits.
\begin{prop}
Given a sequence $\lambda=\{ \lambda_k \}$ of numbers in $[0,1[$, there is a function $P_\lambda$ with zeros at  $\{n+\lambda_m+im, m,n \in \mathbb{Z}\}$, such that $|P_\lambda|(z)\leq e^{c|z|^2}$.
\end{prop}	
{\bf proof}: In the same way as one defines a theta function, we set 
\begin{displaymath}
P_\lambda(z)=\prod_{m=0}^\infty \left(1-e^{-2\pi m}e^{2\pi i\lambda_m} e^{-2\pi i z}\right)\prod_{m=1}^\infty \left (1-e^{-2\pi m}e^{2\pi i\lambda_{-m}} e^{2\pi i z}\right).
\end{displaymath}
Both products are convergent for any given $z\in \mathbb{C}$. In fact, the logarithm of the norm of each of them is bounded by
\begin{displaymath}
\sum_m \log{(1-e^{-2m}e^{2\pi |z|})},
\end{displaymath}
which is a convergent series, because the series $\sum e^{-2m}$ is itself convergent. $P_\lambda$ clearly has zeroes at $n+\lambda_m +im$, as wanted. It remains to be seen that $P_\lambda$ is of order less than 2. First note that $P_\lambda(z+1)=P_\lambda(z)$. It is also true that 
\begin{displaymath}
P_\lambda(z+i)=-e^{-2\pi i z}e^{2\pi}P_{\sigma(\lambda)}(z),
\end{displaymath}
where $\sigma(\lambda)$ is simply a shift left in $\lambda$. This implies
\begin{displaymath}
P_\lambda(z+mi)=(-1)^me^{-2\pi im z}e^{\pi(m^2+m)}P_{\sigma^m(\lambda)}(z)
\end{displaymath}
and from this formula, we can conclude that $P_\lambda$ is of order (less than) 2.
\begin{conj}\label{meuc2}
Given $z_0$ in $\mathbb{C}^*$, $p_k(z_0+\frac{z}{\sqrt{k}})$ has a subsequence converging to a normalization of $P_\lambda$, for some
$\lambda$.
\end{conj}
There is a corresponding conjecture for $q_k$ and $Q_\lambda$ where $Q_\lambda(z)= P_\lambda(-z)$. For the rescaling around zero, the sublimit is a function with zeroes in a spiral. To prove Conjecture \ref{conjectura_principal} using Conjecture \ref{meuc2}, it is enough to establish that $(P_\lambda(z),P_\lambda(-z))$ "satisfies" Theorem \ref{exist_lapiz_holo}. This comes down to checking that $P_\lambda$ has no double zeroes and that $P_\lambda(z)$ and $P_\lambda(-z)$ have no common zeroes (for the first two conditions).

\section{An upper bound for $\eta$ in a symplectic manifold}\label{eta_para_X}
The techniques used in section \ref{eta_para_S2} can actually be generalized to determine upper bounds on how much transversality can be achieved for a "linear system" of dimension $n$ on a $2n$-dimensional symplectic manifold. 

{\bf proof of Theorem \ref{final}}: We will use the same idea as we used to prove this proposition for $S^2$. Namely, the sections $s_0,\cdots,s_n$ define an asymptotically holomorphic map from $X$ to $\mathbb{CP}^n$. The rescaled map subconverges, (in some sense to be specified ahead) as k tends to infinity, to a map of a disc of fixed radius in $\mathbb{C}^n$ to $\mathbb{CP}^n$. If $\eta$ were very close to $1$, the pullback metric from $\mathbb{CP}^n$ would be very close to flat, but this can't be true because its curvature 2-form is the pullback of the Fubini-Study metric.

To make the argument precise, choose a holomorphic structure $J$ on $X$, possibly non integrable, compatible with the symplectic form. In this setting, we know that, for each point $p$ in $X$, there are centered symplectic coordinates, for which $J$ is standard at $p$ and an asymptotically holomorphic sequence of sections $\sigma^k$ of $L^k$ such that 
\begin{displaymath}
||\sigma^k||=e^{-k{|z|^2}},
\end{displaymath}
on a ball of radius $k^{-1/2}$. Although this is explained in \cite{do1}, we will briefly recall it here. Let $\{z_i\}$ be symplectic coordinates near $p$. Choose them so that $J(0)=J_0$ ($J_0$ denotes the standard complex structure on $\mathbb{C}^n$). Then, in coordinates, $\omega=-i\sum dz_jd\bar{z}_j$ and one can choose for the connection 1-form on $L$
\begin{displaymath}
A=\sum z_jd\bar{z}_j-\bar{z}_jdz_j.
\end{displaymath}
The complex structure $J$ allows us to decompose $TB(1)$ into $T^{(0,1)}$ and $T^{(1,0)}$, the eigenspaces for $J$, and to write $d=\partial+\bar{\partial}$. We also have the usual decomposition for $d$, using $J_0$, which we write $d=\partial_0+\bar{\partial}_0$. Let $\mathbf{1}$ denote the unitary section of $L$, satisfying $\nabla \mathbf{1}=A\mathbf{1}$. Then set
\begin{displaymath}
\sigma^k=e^{-k|z|^2}\mathbf{1}^k.
\end{displaymath}  
We have
\begin{displaymath}
\nabla \sigma^k=-2k\sum\bar{z}_jdz_j e^{-k|z|^2}\mathbf{1}^k,
\end{displaymath} 
so that $\bar\partial_0 \sigma^k=0$. But
\begin{displaymath}
|J-J_0|\leq C|z|
\end{displaymath}
and 
\begin{displaymath}
\bar\partial \sigma^k-\bar\partial_0 \sigma^k=i\nabla \sigma^k\circ (J-J_0).
\end{displaymath}
By taking into consideration $|z_j|\leq Ck^{-1/2}$ we get
\begin{displaymath}
|\bar\partial \sigma^k|\leq C.
\end{displaymath}

Write $s^k_l=f^k_l\sigma^k$ and $\tilde{f}^k_l$ for the functions $f^k_l$ on rescaled coordinates, defined over $B(1)$, i.e., $\tilde{f}^k_l=f^k_l\circ \delta_k$, where $\delta_k$ is given in coordinates by
\begin{displaymath}
\delta_k(z)=\frac{z}{\sqrt{k}}.
\end{displaymath}
These functions satisfy
\begin{displaymath}
\eta e^{|z|^2}\leq |\tilde{f}^k_0|^2+\cdots+|\tilde{f}^k_n|^2\leq e^{|z|^2}.
\end{displaymath}
Also, carrying $J$ over to $B(1)$ by using the rescaled coordinates, we get a new complex structure on $B(1)$, which we again call $J$. Let $J_0$ be the standard complex structure on $B(1)$. Because $J$ is $J_0$ at $0$, $J$ and $J_0$ are close. More precisely, they satisfy
\begin{displaymath}
|J-J_0|\leq C\frac{|z|}{\sqrt{k}}, \,\,\,\, z\in B(1).
\end{displaymath}
Again the complex structure $J$ allows us to decompose $TB(1)$ into $T^{(0,1)}$ and $T^{(1,0)}$, the eigenspaces for $J$, and to write $d=\partial+\bar{\partial}$. We also have the usual decomposition for $d$ using $J_0$ which we write $d=\partial_0+\bar{\partial}_0$. We have
\begin{displaymath}
\left| \bar{\partial} \tilde{f^k_l} \right| \leq \frac{C}{\sqrt{k}}.
\end{displaymath}
Next, we show that $\{d\tilde{f}^k_l\}$ is $L^2$ bounded on $B(1/2)$, independently of $k$. Take $\beta$ to be a bump function equal to 1 on $B(1/2)$ and $0$ outside $B(3/4)$. Set $g^k_l=\beta \tilde{f}^k_l$. It is enough to see that 
\begin{displaymath}
\int_{B(1)}|dg^k_l|^2
\end{displaymath}
is bounded. Now,
\begin{displaymath}
\bar\partial g^k_l=\bar\partial\beta \tilde{f}^k_l+\beta \bar\partial\tilde{f}^k_l
\end{displaymath}
and therefore this is uniformly bounded and, as a consequence, $L^2$ bounded. As for 
\begin{displaymath}
\int_{B(1)}|\partial g^k_l|^2,
\end{displaymath}
we want to compare it to 
\begin{displaymath}
\int_{B(1)}|\partial g^k_l|_J^2,
\end{displaymath}
where $|.|_J$ denotes the norm obtained from $J$ and $\omega_0$ (the standard symplectic form on $\mathbb{C}^n$). Given a vector field $v$ on $B(1)$, 
\begin{displaymath}
|v|^2=\omega_0(v,J_0v), \,\,\,\, \,\,\,\,\, |v|_J^2=\omega_0(v,Jv),
\end{displaymath}
so
\begin{displaymath}
\left| |v|^2-|v|_J^2\right|\leq |\omega_0||J-J_0||v|^2\leq C|v|^2
\end{displaymath}
and this implies that the $J$ norm of any tensor which has a bound on its $J_0$ norm, independent of $k$, is bounded independently of $k$ as well. For example, $\omega_0$ and $J_0-J$ have bounded $J$ norms. Conversely, we get
\begin{displaymath}
\left| |v|^2-|v|_J^2\right|\leq |\omega_0|_J |J-J_0|_J |v|_J^2\leq C|v|_{J}^2,
\end{displaymath}
therefore, if a tensor has a bound on its $J$ norm, independent of $k$, it also has a bound (independent of $k$) on its $J_0$ norm. Hence it is enough to prove that $\int_{B(1)}|\partial g^k_l|_J^2$ is bounded. Since $g^k_l$ is zero near the boundary of $B(1)$, we can write
\begin{displaymath}
\int_{B(1)}|\partial g^k_l|_J^2=\int_{B(1)}\partial^*\partial g^k_l \bar{g}^k_l,
\end{displaymath}
where $\partial^*$ means the formal adjoint, with respect to the $J$ metric. The Hodge identities state that $\partial^*=i[\Lambda,\bar\partial]$ and $\bar\partial^*=-i[\Lambda,\partial]$ ($\Lambda$ denotes the dual with respect to the $J$ metric of wedging with $\omega_0$) and 
\begin{displaymath}
\partial^*\partial g^k_l=i\Lambda\bar\partial \partial g^k_l, \,\,\,\, \,\,\,\,\, \bar\partial^*\bar\partial g^k_l=-i\Lambda\partial \bar\partial g^k_l.
\end{displaymath}
The operators $\bar\partial\partial$ and $\partial\bar\partial$ are related by
\begin{displaymath}
\bar\partial\partial+\partial\bar\partial=N\bar{N}+\bar{N} N,
\end{displaymath}
where $N$ is the Nijenhuis tensor for $J$. The upshot is that we can write
\begin{displaymath}
\int_{B(1)}|\partial g^k_l|_J^2=\int_{B(1)}i\Lambda(N\bar{N}+\bar{N} N)|g^k_l|_J^2+\int_{B(1)}\bar\partial^*\bar\partial g^k_l \bar{g}^k_l,
\end{displaymath}
and the second term in this sum is simply $\int_{B(1)}|\bar\partial g^k_l|_J^2$, which we know is bounded. Now $N$ is bounded independently of $k$ and so is $\Lambda$ (because $\omega_0$ is) and we are done.

This implies that $\{g^k_l\}$ is bounded in $W^1(\mathbb{C}^n)$ and so, by Rellich's lemma, it has a subsequence, converging to $f_l$ in $W^0=L^2$. Now by Riesz theorem, we can extract from $\{g^k_l\}$ a subsequence converging pointwise a.e. in $B(1)$ and $\{\tilde{f}^k_l\}$ has a pointwise convergent subsequence a.e. in $B(1/2)$ to $f_l$. These functions satisfy
\begin{displaymath}
\eta e^{|z|^2}\leq |{f}_0|^2+\cdots+|{f}_n|^2\leq e^{|z|^2}, \,\, a.e. \,.
\end{displaymath}

The next step is to show that the functions $f_l$ are holomorphic. To do this, note first that $\{\tilde{f}^k_l\}$ has a weakly convergent subsequence in $W^1$, simply because it is a bounded sequence in this space. Given any $L^2$ 1-form $a$ in $B(1)$,
\begin{displaymath}
\int_{B(1/2)}\langle d\tilde{f}^k_l-df_l,a \rangle \rightarrow 0.
\end{displaymath}
In fact, if $a\in L^2(T_0^{(0,1)}B(1))$
\begin{displaymath}
\int_{B(1/2)}\langle\bar\partial_0\tilde{f}^k_l,a \rangle \rightarrow \int_{B(1/2)}\langle\bar\partial_0 f_l,a \rangle 
\end{displaymath}
because $T_0^{(0,1)}$ and $T_0^{(1,0)}$ (the eigenspaces of $J_0$) are pointwise orthogonal. But
\begin{displaymath}
\bar\partial_0\tilde{f}^k_l=(\bar\partial_0-\bar\partial) \tilde{f}^k_l+\bar\partial\tilde{f}^k_l,
\end{displaymath}
so, over $B(1/2)$,
\begin{displaymath}
||\bar\partial_0\tilde{f}^k_l||_{L^2}\leq||d\tilde{f}^k_l\circ(J-J_0)||_{L^2}+||\bar\partial\tilde{f}^k_l||_{L^2}.
\end{displaymath}
The first term is bounded by
\begin{displaymath}
\frac{C}{\sqrt{k}} \left( \int_{B(1/2)}|d\tilde{f}^k_l|^2\right)^2\leq \frac{C}{\sqrt{k}}.
\end{displaymath}
We know that the second term is also bounded by $C/\sqrt{k}$. We can conclude that 
\begin{displaymath}
||\bar\partial_0\tilde{f}^k_l||_{L^2}\rightarrow 0,
\end{displaymath}
so
\begin{displaymath}
(\bar\partial_0\tilde{f}^k_l,a)_{L^2}\rightarrow 0.
\end{displaymath}
This in turn proves that 
\begin{displaymath}
(\bar\partial_0{f}_l,a)_{L^2}=0, \,\,\,\, \forall a\in L^2(T_0^{(0,1)}B(1))
\end{displaymath}
and, in fact, for all $a\in L^2(T^*B(1))$, so $f_l$ is holomorphic. In particular, it is continuous and the inequality 
\begin{displaymath}
\eta e^{|z|^2}\leq |{f}_0|^2+\cdots+|{f}_n|^2\leq e^{|z|^2}, 
\end{displaymath}
holds at every point. We conclude that the existence of the $n$ sections $s^k_0,\cdots,s^k_n$ implies the existence of $n$ holomorphic functions $f_0,\cdots,f_n:B(1/2)\rightarrow \mathbb{C}$ satisfying
\begin{displaymath}
\eta e^{|z|^2}\leq |{f}_0|^2+\cdots+|{f}_n|^2\leq e^{|z|^2}.
\end{displaymath}
Next, we show that the existence of such functions implies the equality
\begin{equation}\label{d_barra_d=d_barra_d}
\triangle \log(|f_0|^2+\cdots+|f_n|^2)=\triangle \log \det Hess(|f_0|^2+\cdots+|f_n|^2),
\end{equation}
at those points where the map $F=[f_0:\cdots:f_n]$ has injective differential, so that the pullback of the Fubini-Study metric is a metric. 
Away from its branch points, ${F}$ pulls back the Fubini-Study metric on $\mathbb{CP}^n$ to a metric on $B(c)$, whose curvature 2-form is $F^*\omega_{FS}$. On the other hand, this metric is the metric associated to the 2-form $F^*\omega_{FS}$. On $\{[\mathbf{z}_0:\cdots :\mathbf{z}_n] :\mathbf{z}_0\ne 0\}\subset \mathbb{CP}^n$ there are inhomogeneous coordinates 
\begin{displaymath}
(z_1=\frac{\mathbf{z}_1}{\mathbf{z}_0},\cdots ,z_n=\frac{\mathbf{z}_n}{\mathbf{z}_0}), 
\end{displaymath}
and the Fubini-Study metric can be written as
\begin{displaymath}
\bar{\partial}\partial \log(1+|z_1|^2+\cdots+|z_n|^2).
\end{displaymath}
The map $F$ being holomorphic, we conclude that $F^*\omega_{FS}$ is
\begin{displaymath}
\bar{\partial}\partial \log(1+|z_1|^2+\cdots+|z_n|^2)\circ F.
\end{displaymath}
Now 
\begin{displaymath}
\log(1+|z_1|^2+\cdots+|z_n|^2)\circ F=\log(|f_0|^2+|f_1|^2+\cdots+|f_n|^2)-\log(|f_0|^2),
\end{displaymath}
so
\begin{displaymath}
F^*\omega_{FS}=\bar{\partial}\partial \log(|f_0|^2+|f_1|^2+\cdots+|f_n|^2).
\end{displaymath}
Next we prove a lemma:
\begin{lemma}
Consider the metric on $B(c)$, the ball of radius $c$ around the origin in $\mathbb{C}^n$, associated with the 2-form $\bar{\partial}\partial u$, for some plurisubharmonic function $u$ on $B(c)$. Then, its curvature 2-form is 
\begin{displaymath}
\bar{\partial}\partial \log \det \text{Hess} (u)
\end{displaymath}
where $Hess(u)$ denotes the matrix $(u_{z_i\bar{z_j}})$.
\end{lemma}
{\bf proof}: On the bundle $TB(c)$ we have a trivialization given by 
\begin{displaymath}
(\partial/\partial z_1,\cdots ,\partial/\partial z_n).
\end{displaymath}
With respect to this trivialization the matrix representation of the metric is $h=Hess(u)$. There is a unique connection on $TB(c)$ which is compatible with the metric associated with $\bar{\partial}\partial u$ and the complex structure. Its matrix representation, with respect to the mentioned trivialization of the bundle, is $\theta=\partial h h^{-1}$. The curvature matrix representation is $d\theta+\theta \wedge \theta$. The curvature 2-form is the trace of this two by two matrix of 2-forms. Since $\theta \wedge \theta$ is traceless, the curvature is $tr d\theta=d tr \theta$. But we know that
\begin{displaymath}
tr(\partial h h^{-1})=\partial \log(\det h).
\end{displaymath}
The curvature becomes $d\partial \log(\det h)=\bar\partial \partial \log(\det h)$, as we wished to prove. 

Going back to our problem, we can apply this lemma to $u= \log(|f_0|^2+\cdots+|f_n|^2)$, to conclude that the curvature of the pullback by $F$ of the Fubini-Study metric is 
\begin{displaymath}
\bar\partial \partial \log(\det Hess(\log(|f_0|^2+\cdots+|f_n|^2)),
\end{displaymath}
away from the branch points of $F$ and the zeroes of $f_0$. We therefore conclude that
\begin{displaymath}
\bar\partial \partial \log(\det Hess(\log(|{f}_0|^2+\cdots+|{f}_n|^2))=\bar\partial \partial \log(|f_0|^2+\cdots+|f_n|^2),
\end{displaymath}
where $F$ is an embedding. 
If this set of points is dense in $B(1)$, we are done. In fact, then, equality (\ref{d_barra_d=d_barra_d}) holds true at all points. Assume that $\eta$ is very close to 1, then 
\begin{displaymath}
\log(|f_0|^2+\cdots+|f_n|^2)\simeq |z|^2,
\end{displaymath}
so that 
\begin{displaymath}
\triangle \log(|f_0|^2+\cdots+|f_n|^2)\simeq n,
\end{displaymath}
and
\begin{displaymath}
\log \det Hess \log(|f_0|^2+\cdots+|f_n|^2)\simeq 0.
\end{displaymath}
This gives a contradiction. At first, it seems that the upper bound $\eta_0$ obtained for $\eta$ in this way depends on the particular $f_0,\cdots,f_n$ and therefore on the specific sequences $s^k_0,\cdots,s^k_n$, but this is not so. Using Lemma \ref{delta_lambda} (and a slight variation of it applicable to $Hess$ instead of $\triangle$), the only thing we need to note here is that there is a bound (independent of the sequences) on
\begin{displaymath}
d\triangle  \log(|f_0|^2+\cdots+|f_n|^2),
\end{displaymath}
\begin{displaymath}
d Hess \log(|f_0|^2+\cdots+|f_n|^2),
\end{displaymath}
and on
\begin{displaymath}
d\triangle \log \det Hess \log(|f_0|^2+\cdots+|f_n|^2).
\end{displaymath}
This is so because, since $f_0,\cdots,f_n$ are holomorphic, the Cauchy formula and the bounds $|f_i|^2\leq e^{c^2}$, give bounds on all the derivatives of the $f_i$ depending only on $c$, and, choosing $c$ appropriately small, the above quantities have bounds depending only on $n$, just as for $\mathbb{CP}^1$. 

On a maybe smaller ball, the differential of $F$ is injective exactly where the determinant of the derivative of $(f_1/f_0,\cdots,f_n/f_0)$ is non zero. But this determinant is a holomorphic function, so, if it is not identically zero, its zero set will have zero measure. In the case where the determinant is identically zero, we make use of the following: given any $\epsilon$, there are holomophic functions $\alpha_1,\cdots, \alpha_n$, each with norm smaller than $\epsilon$, such that 
\begin{displaymath}
\det d\left( (\frac{f_1}{f_0},\cdots,\frac{f_n}{f_0})+(\alpha_1,\cdots,\alpha_n)\right)
\end{displaymath}
is not identically zero. Now set $g_0=f_0$ and $g_i=f_i+f_0\alpha_i$, for $i$ in $\{1,\cdots n\}$. Then, by choosing $\epsilon$ small, we can ensure that
\begin{displaymath}
(2\eta-1) e^{|z|^2}\leq |{g}_0|^2+\cdots+|{g}_n|^2\leq (2-\eta)e^{|z|^2}.
\end{displaymath}
Since the $g_i$ are holomorphic, by the same reasoning, we get an upper bound on $\eta$ independent of the $f_i$'s.

\section{Further remarks}

 Let $(p_k,q_k)$ be polynomials of degree smaller than $k$ in one complex variable $z$. Set
\begin{displaymath}
\rho_k(z)=\frac{|p_k|^2(z)+|q_k|^2(z)}{(1+|z|^2)^{2k}}.
\end{displaymath} 
That is, $\rho_k$ is the sum of the squares of the norms of $p_k$ and $q_k$, seen as sections of $\mathcal{O}(k)$. If these polynomials satisfy the conditions in Proposition \ref{principal} (using inhomogeneous coordinates in $\mathbb{CP}^1$), the function $\rho_k$ is bounded above by 1 and below by $\eta$. We have seen that the function $v_k=-\log(\rho_k)$ satisfies
\begin{displaymath}
\triangle v_k=k-K_ke^{2v_k},
\end{displaymath}
where $K_k$ is a positive function vanishing at the branch points of $p_k/q_k$ and completely determined by those branch points. The PDE,
\begin{displaymath}
\triangle v=c-Ke^{2v},
\end{displaymath}
has been extensively studied when $c=1$ (see \cite{kw} and \cite{cy}). It is the equation for prescribing curvature on $S^2$. There are variational methods adapted to this equation when $c$ is a small constant but these do not apply to our case. Proposition \ref{principal} implies:
\begin{prop}
There is a constant $C$, such that for every $k$ big enough, there is a set of  $2k-2$ points such that the equation
\begin{displaymath}
\triangle v_k=k-K_ke^{2v_k}
\end{displaymath}
has a solution $v_k$ with $||v_k||_{\infty}\leq C$.
\end{prop}
It would be interesting to see if one could recover this result from a PDE theory point of view. This would give an alternative proof of Donaldson's main theorem in \cite{do2} for $S^2$ and maybe a way to calculate $\eta$, by calculating $C$.

Consider the following functionals on the space of $\mathcal{C}^{\infty}$ maps from $S^2$ to $S^2$:
\begin{displaymath}
E_1(F)=\int\left|\frac{F^*\omega}{\omega}\right|^2,
\end{displaymath}
and
\begin{displaymath}
E_2(F)=\int\left|\frac{F^*\omega}{\omega}\right|^4.
\end{displaymath}
It is not hard to see that, the restriction of these two functionals to the space of rational maps of degree $k$ from $S^2$ to itself, is proper and therefore has a minimum. We conjecture that the minimizers of such functionals arise as quotients of polynomials satisfying the property in Proposition \ref{principal}. Similarly, if one considers the much studied functional,
\begin{displaymath}
\int\left|dF\right|^2+\epsilon\int\left|dF\right|^4,
\end{displaymath}
restricted to the space of rational maps of degree $k$, it is natural to conjecture that its minimum also arises in this way.


\begin{thebibliography}{pla1}

\bibitem[Au1]{a1} D. Auroux, {\em Asymptotically holomorphic families of symplectic submanifolds }, Geom. Funct. Anal.  {\bf 7} (1997), 971-995.

\bibitem [BSS]{s} S. B\"ocherer, P. Sarnak, R. Schulze-Pillot, {\em Arithmetic and
Equidistribution of measures on the sphere}, Commun. Math. Phys.
{\bf 242} (2003), 67-80.

\bibitem[CY]{cy} A. Chang, P. Yang, {\em Prescribing Gaussian curvature on $S^2$} , Acta. Math., {\bf159} (1987), 214-259. 

\bibitem[Do2]{doa} S.K. Donaldson, {\em Lefschetz Fibrations in
Symplectic Geometry}, Documenta Mathematica Extra Volume ICM (1998)
, II, 309-314.

\bibitem[Do1]{do1} S.K. Donaldson, {\em Symplectic submanifolds and
almost-complex geometry}, J. Differential Geometry {\bf 44} (1996),
666-705.

\bibitem[Do3]{do2} S.K. Donaldson, {\em Lefschetz pencils on symplectic
manifolds}, J. Differential Geom. {\bf 53} (1999), 205-236.

\bibitem [DS]{sd} S.K. Donaldson, I. Smith, {\em Lefschetz pencils and
the canonical class for symplectic 4-manifolds}, Topology {\bf 42}
(2003), 743-785.

\bibitem [KW]{kw} J. Kazdan, F. Warner, {\em Curvature functions for compact 2-manifolds}, Ann. of Math. (2) {\bf 99} (1974), 14Ð47. MR 49:7949 

\bibitem[LPS]{lps}  A. Lubotsky, R. Phillips, P. Sarnak, {\em Hecke operators and distributing points on $S^2$ I and II}, Comm. Pure and Appl. Math., {\bf 34}, (1986), 149-186, and No. XL, (1987), 401-420. 

\bibitem[RSZ]{rsz} E.A. Rakhmanov, E.B. Saff, Y.M. Zhou, {\em Minimal
discrete energy on the sphere}, Math. Research Letters {\bf 1}
(1994), 647-662.

\end{thebibliography}
\end{document}